\documentclass[10pt,a4paper]{article} 

\usepackage[latin1]{inputenc}  
\usepackage[babel=true]{csquotes}
\usepackage{amsmath,amscd,amssymb,latexsym,graphicx,color,mathrsfs} 
 
 \usepackage[all]{xy}
\newtheorem{theorem}{Theorem}[section]

\newtheorem{corollary}[theorem]{Corollary}

\newtheorem{definition}[theorem]{Definition}

\numberwithin{equation}{section} 
 
\newcommand{\cqfd}{\hfill{\small $\Box$}} 
 \newenvironment{proof}[1][]{{\bf Proof #1 : }}{\hfill \cqfd} 
 
\reversemarginpar

\newcommand{\cn}{\mathrm{cn}} 
\newcommand{\pcn}{\mathrm{pcn}}

\def\intx{\overset{\:\circ}{X}}
\newcommand\tgt[1]{{}^{T}\kern-1pt #1}
\newcommand\adi[1]{{}^{ad}\kern-1pt #1}

  \def\CC{{\mathbb{C}}}

  \def\RR{{\mathbb{R}}}

 \def\ZZ{{\mathbb{Z}}}

  \def\bR{{\mathbf{R}}}

\def\cA{{\mathcal{A}}}  
  \def\cF{{\mathcal{F}}}
 \def\cH{{\mathcal{H}}} 
 \def\cK{{\mathcal{K}}} 
  
\def\cP{{\mathcal{P}}}  
\def\cS{{\mathcal{S}}}

\title{Fredholm anomalies on manifolds with corners of low codimensions and conormal corner cycles \footnote{Both authors were partially supported by the ANR OpART, ANR-23-CE40-0016.} }

\author{P. Carrillo Rouse and J.M. Lescure} 
 
\date{\today} 
 
\begin{document}

\maketitle 

\begin{center}
{\bf Abstract}
\end{center}

Given a connected manifold with corners $X$ of any codimension there is a very basic and computable homology theory called conormal homology defined in terms of faces and orientations of their conormal bundles, and whose cycles correspond geometrically to corner's cycles, these conormal homology groups are denoted by $H^{cn}_*(X)$. Using our previous works we define an index morphism
$$K^0(^bT^*X)\stackrel{Ind_{ev,cn}^X}{\longrightarrow}H_{ev}^{cn}(X)$$
for $X$ a manifold with corners of codimension less or equal to three and called here the even conormal index morphism. In the case that $X$ is compact and connected and $D$ is an elliptic $b-$pseudodifferential operator in the associated $b-$calculus of $X$ we know, by our previous works and other authors works, that, up to adding an identity operator, $D$ can be perturbed (with a regularizing operator in the calculus) to a Fredholm operator iff $Ind_{ev,cn}^X([\sigma_D])$ (where $[\sigma_D]\in K^0(^bT^*X)$ is the principal symbol class) vanishes in the even conormal homology group $H_{ev}^{cn}(X)$.
The main result of this paper is the explicit computation of the even and odd conormal index morphisms $Ind_{ev/odd,cn}^X(\sigma)\in H_{ev/odd}^{cn}(X)$ for $X$ a manifold with corners of codimension less or equal to three. The coefficients of the conormal corner cycles $Ind_{ev/odd,cn}^X(\sigma)$ are given in terms of some suspended Atiyah-Singer indices of the maximal codimension faces of $X$ and in terms of some suspended Atiyah-Patodi-Singer indices of the non-maximal codimension faces of $X$. As a corollary we give a complete caracterization to the obstruction of the Fredholm perturbation property for closed manifolds with corners of codimension less or equal to three in terms of the above mentioned indices of the faces, this allows us as well to give such a characterization in terms of the respective topological indices.

\tableofcontents

\section{Introduction}
We continue our study on obstructions on Fredholm boundary conditions on manifolds with corners initiated in \cite{CarLes} and in \cite{CRLV}.

Given a connected manifold with corners $X$ of any codimension there is a very basic and computable homology theory called conormal homology defined in terms of faces and orientations of their conormal bundles, and whose cycles correspond geometrically to corner's cycles. These conormal homology groups are denoted by $H^{cn}_*(X)$ in their integer coefficients version. The main theorem in \cite{CarLes} gives a natural and explicit isomorphism
\begin{equation}\label{Tisointro}
\xymatrix{
K_*(\cK_b(X))\ar[r]^-T_-\cong & H^{pcn}_*(X)
}
\end{equation}
between the $K-$theory group (even or odd) of the algebra $\cK_b(X)$ of $b$-compact operators for $X$ and  the periodic (even or odd) conormal homology group with integer coefficients for the case of $X$ a manifold with corners of codimension less or equal to three (or products of such manifolds with corners). Later on, together with Mario Velasquez in \cite{CRLV}, we gave, for $X$ of arbitrary codimension, an explicit and natural morphism
\begin{equation}\label{Trationalintro}
\xymatrix{
K_*(\cK_b(X))\ar[r]^-T& H^{pcn}_*(X;\mathbb{Q})
}
\end{equation}
and we proved it is a rational isomorphism. We will recall below why the techniques used in the couple of references mentioned above are very different and we will in particular explain why this paper deals only with lower codimensions. But first, let us recall that the main motivation for computing the K-theory groups $K_*(\cK_b(X))$ is because they are the receptacle of higher indices of $b$-elliptic operators and that these higher indices measure the obstruction of such operators to be Fredholm up to some perturbation. Let us recall this in more detail since it is at the heart of the content of this article.

On a smooth compact manifold, ellipticity of (classical) pseudodifferential operators is equivalent to Fredholmness, and the vanishing of the Fredholm index of an elliptic  pseudodifferential operator is equivalent to its invertibility after perturbation by a regularizing operator. In the case of a smooth manifold with boundary, not every elliptic operator is Fredholm and it is known since Atiyah and Bott that there exist obstructions to  the existence of local boundary conditions in order to upgrade an elliptic operator into a Fredholm boundary value problem. Nonetheless, if one moves to non-local boundary conditions, obstructions disappear: for instance, not every elliptic (b-operator)   pseudodifferential operator is Fredholm but it can be perturbed with a regularizing operator to become Fredholm. This non trivial fact, which goes back to Atiyah, Patodi and Singer \cite{APS1}, can also be obtained from the vanishing of a boundary analytic index 
(see \cite{Mel-Pia1997-1,Mel-Pia1997-2,MontNis}, and below). In fact, in this case the boundary analytic index takes values in the $K_0$-theory group of the algebra of regularizing operators and this K-theory group is easily  seen to vanish. It is known that obstructions to the existence of perturbations of elliptic operators into Fredholm ones reappear in the case of manifolds with corners of arbitrary codimension (\cite{Bunke,NSS2010}) (this includes for instance many useful domains in Euclidean spaces). 


Using K-theoretical tools for solving index problems was the main asset in the series of papers by Atiyah-Singer (\cite{AS,AS3}).  In the case of manifolds with boundary, $K$-theory is still useful to understand the vanishing of the obstruction to the existence of pertubations of elliptic operators into Fredholm ones (even if $K$-theory is not essential in the computation of analytical indices \cite{APS1}), and a fortiori to understand this obstruction in the case of families of manifolds with boundary (\cite{Mel-Pia1997-1,Mel-Pia1997-2,Mel-Roc2006}).  
For manifolds with corners, Bunke \cite{Bunke} has delivered  for Dirac type operators a complete study of the obstruction, which lives in the homology of a complex associated with the faces of the manifold. As  observed in $\cite{CarLes}$, this homology also appears as the $E^2$-term of the spectral sequence computing the $K$-group that contains the obstruction to Fredholmness for general elliptic $b$-pseudodifferential operators. Nazaikinskii, Savin and Sternin also use $K$-theory to express the obstruction for their pseudodifferential calculus on manifolds with corners and stratified spaces \cite{NSS2,NSS2010}.

Let us briefly recall the framework in which we are going to work.  The algebra of pseudodifferential operators $\Psi_b^*(X)$ associated to any manifold with corners $X$ is defined in \cite{MelPia}: it generalizes the case of manifolds with boundary treated in \cite{Mel} (see also \cite[Section 18.3]{Horm-3}).  The elements in this algebra are called $b-$pseudodifferential operators\footnote{To simplify we discuss only the case of scalar operators, the passage to operators acting on sections of vector bundles is done in the classic way.}, the subscript $b$ identifies these operators as obtained by microlocalization of the Lie algebra of $C^\infty$ vector fields on $X$ tangent to the boundary. This Lie algebra of vector fields can be explicitly obtained as sections of the so called $b$-tangent bundle $^bTX$ (compressed tangent bundle that we will recall below).  The b-pseudodifferential calculus has the classic and expected properties. In particular there is a principal symbol map
$$\sigma_b:\Psi_b^m(X)\to S^{[m]}(^bT^*X).$$  
Ellipticity has the usual meaning, namely invertibility of the principal symbol. Moreover (discussion below and Theorem 2.15 in \cite{MelPia}), an operator is elliptic if and only\footnote{Notice that this remark implies that to an elliptic $b$-pseudodifferential operator one can associate an "index" in the algebraic K-theory group $K_0(\Psi_b^{-\infty}(X))$ (classic construction of quasi-inverses).} if it has a quasi-inverse modulo $\Psi_b^{-\infty}(X)$. Now,  $\Psi_b^{-\infty}(X)$ also contains non compact operators and compacity is there characterized by the vanishing of a suitable indicial map (p.8 ref.cit.). Elliptic $b$-pseudodifferential operators being invertible modulo compact operators -and hence Fredholm\footnote{see p.8 in \cite{MelPia} for a characterization of Fredholm operators in terms of an indicial map or \cite{Loya} 
theorem 2.3 for the proof of Fully ellipticity iff Fredholm}-, are usually said to be  {\sl fully} elliptic. 

The norm  closure $\cK_b(X)$  of  $\Psi_b^{-\infty}(X)$ into the bounded operators on $L^2_b(X)$ fits in the short exact sequence of $C^*$-algebras:
\begin{equation}\label{Introbcompact}
\xymatrix{
0\ar[r]&\cK(X)\ar[r]^-{i_0}&\cK_b(X)\ar[r]^-{r}&\cK_b(\partial X)\ar[r]&0,}
\end{equation}
where $\cK(X)$ denotes the algebra of compact operators on $L^2(X)$ and where $\cK_b(\partial X)$ can be considered to be a quotient algebra for the purposes of this introduction. 

In order to understand how the above sequence enters into the study of Fredholm Perturbation properties we need to settle some definitions.

{\bf Analytic Index morphism:}  Every $b$-pseudodifferential operator has  a principal symbol class in $K^0_{top}(^bT^*X)$, and possesses an  interior parametrix that brings a class in $K_0(\cK_b(X))$ called the {\sl analytical index class}. Both classes are related as follows. Consider the short exact sequence 
\begin{equation}\label{IntrobKses}
\xymatrix{
0\ar[r]&\cK_b(X)\ar[r]&\overline{\Psi_b^0(X)}\ar[r]^-{\sigma_b}&C(^bS^*X)\ar[r]&0.
}
\end{equation}
After applying the $K$-functor, it gives rise to the boundary morphism:
$K_1(C(^bS^*X))\to K_0(\cK_b(X))$ that can be factorized canonically into a  morphism
\begin{equation}\label{intro:ana-ind}
\xymatrix{
K^0_{top}(^bT^*X)\ar[r]^-{Ind^a_X}&K_0(\cK_b(X))
}
\end{equation}
called {\it the Analytic Index morphism of $X$}, which is the one that maps the principal symbol class to the analytical index class of a given elliptic $b$-operator. Now, when composed with the isomorphism (\ref{Tisointro}) , we obtain a morphism

\begin{equation}\label{intro:evencn-ind}
\xymatrix{
K^0_{top}(^bT^*X)\ar[rr]^-{Ind^X_{ev,cn}}&&H^{cn}_{ev}(X)
}
\end{equation}
called here {\bf The even conormal index morphism of $X$}. This last index morphism is  the main subject of study in this paper as explained below.

To express how \eqref{Introbcompact} and the previous index maps relate to the Fredholm obstruction we introduce the following vocabulary:
\begin{definition} Let $D\in \Psi_b^m(X)$ be elliptic. We say that $D$ satisfies: 
\begin{itemize}
\item  the {\it Fredholm Perturbation Property} $(\cF\cP)$ if there is   $R\in \Psi_b^{-\infty}(X)$ such that $D+R$ is  fully elliptic.  
\item the  {\it stably  Fredholm Perturbation Property} $(\cS\cF\cP)$ if $D\oplus 1_H$ satisfies $(\cF\cP)$ for some identity operator $1_H$.
\end{itemize}
\end{definition}
 The following theorem is a corollary of theorem 2.22 in \cite{CRLV} and of the main computation in \cite{CarLes}, that is, the fact that $T$ in (\ref{Tisointro}) is an isomorphism. In its abstract version it was proved first in a weaker form in \cite{NSS2} (see \cite{CRLV} Theorem 2.22 for more details).

\begin{theorem}\label{thmintroconormalobstruction}
Let $D$ be an elliptic $b$-pseudodifferential operator on a compact manifold with corners $X$ of codimension lower or equal to three. Then  $D$ satisfies $(\cS\cF\cP)$ if and only if the even conormal index 
$
 Ind_X^{cn,ev}([\sigma_b(D)])=0$ in  $H_{ev}^{cn}(X)$.  \\ In particular, if $D$ satisfies $(\cF\cP)$ then its even conormal index vanishes.
\end{theorem}

The main result of this paper is the explicit computation of $Ind_X^{cn,ev}(\sigma)$ as an element in $H_{ev}^{cn}(X)$ together with the corollaries explained below. Now, conormal homology is defined by a very elementary chain complex, indeed, for the purpose of this introduction, if $X$ is a manifold with corners of codimension $d$ and $p\in \{0,...,d\}$, the $p$-conormal chains $C_p(X)$ is essentially (see section \ref{subsectionconormalhomology} for precisions) the $\mathbb{Z}-$module generated by $F_p(X)$=faces of codimension $p$ and the differential  $\delta_p : C_p(X)\to C_{p-1}(X)$ looks like (complete definition in \ref{subsectionconormalhomology})
\begin{equation}\label{diffpcn1}
   \delta_p(f) = \sum_{\substack{g\in F_{p-1},  \\ f\subset\overline{g}}} \sigma(f,g)\cdot g
\end{equation}
where $\sigma(f,g)=\pm 1$ depending on the some orientations given by the conormal bundles. The {\it conormal homology} of $X$, denoted by $H^{\cn}_*(X)$ is defined to be the homology of $(C_*(X),\delta_*)$. Even and odd groups are called  {\it periodic conormal homology}:
 \begin{equation}
 H^{cn}_{ev}(X)=\oplus_{p\ge 0} H^{\cn}_{2p}(X) \text{ and }H^{cn}_{odd}(X)=\oplus_{p\ge 0} H^{\cn}_{2p+1}(X).
\end{equation}
In particular, an element $\omega\in H_{2p}^{cn}(X)$ is represented by a corner cycle $\sum_{Y\in F_{2p}}a_Y\cdot Y\in Ker\, \delta_{2p}$, and computing $Ind_X^{cn,ev}(\sigma)\in H_{ev}^{cn}(X)$ means to us to be able to describe the coefficients $a_Y\in \mathbb{Z}$, this is what we do and explain below.

The last theorem justifies our interest for computing the even conormal indices. In fact, we will actually compute all the even and odd conormal indices up to codimension three closed or not\footnote{These indices make sense for not necessarilly closed mwc as seen below and coincide with their analytically defined versions when the manifolds are closed.}, we do this as an independently interesting algebraic topology and K-theoretical exercise that will apply, as seen below, to get explicit Fredholm obstructions results together with topological interpretations. Elsewhere we will see that the odd indices are also very interesting to study.

To state the explicit computations we need to introduce some terminology and in particular some indices that will appear on the coefficients of the conormal cycles we are going to compute.

First of all, we fix for the entire paper a generator $B$  of $K_1(C_0(\mathbb{R}))$ and let 
$$\beta_p:=B^p =B \otimes_{\CC}\cdots\otimes_{\CC}B	 \in K_p(C_0(\RR^p)).$$
This choice is important since it fits with the explicit computation of the isomorphism (\ref{Tisointro}) done in \cite{CarLes} and recalled below in \ref{subsectionKvsconormal}.

Also, for keeping the notation short, for a Lie groupoid $\mathscr{G}$ we let $K^*(\mathscr{G})$ denote the topological $K$-theory group $K_*(C^*(\mathscr{G}))$ of the Lie groupoid $C^*-$algebra (max or reduced, all the groupoids used in this paper are amenable so these algebras coincide).

Now, let $Y$ be a smooth manifold. For $p\in \mathbb{N}$ we define the suspended (or $p$-suspended if needed) Atiyah-Singer index morphism associated to $Y$ to be the morphism ($*=0$ if $p$ even or $*=1$ if $p$ odd)
\begin{equation}
\xymatrix{
K^*(T^*Y\times \mathbb{R}_p)\ar[rr]^-{Ind_{AS,p}^Y}&& \mathbb{Z}
}
\end{equation} 
where $\mathbb{R}_p$ stands above for the additive group $(\mathbb{R}^p,+)$ and defined 
as the morphism fitting in the following commutative diagram
\begin{equation}
\xymatrix{
K^*(T^*Y\times \mathbb{R}_p)\ar[d]_-{\beta_p^{-1}}\ar[rr]^-{Ind_{AS,p}^Y}&& \mathbb{Z}\ar[d]^-=\\
K^0(T^*Y)\ar[rr]_-{Ind_{AS}^Y}&& \mathbb{Z}
}
\end{equation}
where $Ind_{AS}^Y$ is the classic Atiyah-Singer index morphism\footnote{This index can be defined for any $Y$, closed or not, the analytic interpretation of it is not the same when $Y$ is not closed.} for $Y$ and where $\beta_p:K^0(T^*Y)\to K^*(T^*Y\times \mathbb{R}_p)$ is the composition of the Bott isomorphism $K^0(T^*Y)\to K^*(T^*Y\times \mathbb{R}^p)$ given by external multiplication with the Bott element $\beta_p$ described above, followed by the K-theory isomorphism induced by the isomorphism of algebras 
$$C^*(T^*Y\times \mathbb{R}^p)\cong C^*(T^*Y\times \mathbb{R}_p)$$
induced by the Fourier isomorphism in the direction $\mathbb{R}^n$.

Now, let $Z$ be a manifold with corners. As recalled below (see sections \ref{subsectionFredindexgeo} and \ref{Fredsubsection} for more details), there is an index morphism
\begin{equation}
\xymatrix{
K^0(^bT_{nc}Z)\ar[rr]^-{Ind_{Fred}^Z}&& \mathbb{Z}
}
\end{equation}
computing the Fredholm index of a fully elliptic (equivalent to Fredholm) operator in the case the manifold with corners $Z$ is closed. Indeed, as explained below (see sections \ref{subsectionFredindexgeo} and \ref{Fredsubsection} for more details) any such operator has a noncommutative symbol defining a class in the $K_0$-group of the $C^*$-algebra of a deformation groupoid $^bT_{nc}Z$ called the noncommutative tangent groupoid, see as well \cite{DLN} and \cite{DLR} for further motivation on the terminology. Now, for $p\in \mathbb{N}$ we define the suspended (or $p$-suspended if needed) Fredholm index morphism associated to $Z$ to be the morphism ($*=0$ if $p$ even or $*=1$ if $p$ odd)
fitting the commutative diagram
\begin{equation}
\xymatrix{
K^*(^bT_{nc}Z\times \mathbb{R}_p)\ar[d]_-{\beta_p^{-1}}\ar[rr]^-{Ind_{Fred,p}^Z}&& \mathbb{Z}\ar[d]^-=\\\
K^0(^bT_{nc}Z)\ar[rr]_-{Ind_{Fred}^Z}&& \mathbb{Z}
}
\end{equation}
where $Ind_{Fred}^Z$ is the Fredholm index morphism\footnote{Again, these indices can be defined for any mwc $Z$, closed or not. As recalled below, they coincide with their analytically defined versions when $Z$ is closed. Still, as for the case of all the indices in this paper, it is very important to us that these indices can be constructed for not necessarilly closed mwc. In \cite{CLM}, for the case of manifolds with boundaries, this index morphism was called the APS index morphism associated to such a manifold, the name was inspired by the foundational work of Atiyah, Patodi and Singer \cite{APS1} and its sequels.} of the manifold with corners $Z$ and where $\beta_p$ is the composition of the Bott isomorphism $K^0(^bT_{nc}Z)\to K^*(^bT_{nc}Z\times \mathbb{R}^p)$ given by external multiplication with the Bott element $\beta_p$ described above, followed by the isomorphism induced by the Fourier isomorphism 
$$C^*(^bT_{nc}Z\times \mathbb{R}^p)\cong C^*(^bT_{nc}Z\times \mathbb{R}_p).$$

We are ready to state the main theorem of this paper, it resumes the content of the entire section \ref{sectionevenoddconormalindex} where the computations of the even/odd indices are done codimension by codimension. In the statement there are some important elements associated to a symbol class $\sigma\in K^*(^bT^*X)$, for example there are obvious restricted symbols $\sigma\mapsto \sigma_Y$ where $Y$ is a maximal codimension face of $X$. But there are also some associated (suspended) non commutative symbols
$$\sigma_{0,Z}^{nc}\in K^*(T_{nc}\overline{Z}\times \mathbb{R}_{codim(Z,X)}),$$
where $Z$ is a (open) face of not maximal codimension $codim(Z,X)$ (codimension of the face $Z$ in the manifold with corners $X$) and hence $\overline{Z}$ is a manifold with corners of codimension $codim(Z,X)$.
These noncommutative symbols are not at all arbitrary and their construction is essential for the statement and for the proof of the computations below.

For precise notations and the precise definition of these noncommutative and commutative symbols see section \ref{sectionevenoddconormalindex}. The theorem states as follows:

\begin{theorem}
Let $X$ be a manifold with corners of codimension less or equal to three. We have
\begin{itemize}
\item {\bf Odd conormal index in codim 1\footnote{As discussed above, the even conormal index for a $1$-codimensional manifold with corners is zero.}}: If $codim(X)=1$ and $\sigma\in K^1(^bT^*X)$ we have that
\begin{equation}
\sum_{Y\in F_1(X)}Ind_{AS,1}^{Y}(\sigma_Y)\cdot Y \in Ker\,\delta_1
\end{equation}
and 
\begin{equation}
Ind_{cn,odd}^X(\sigma)=\sum_{Y\in F_1(X)}Ind_{AS,1}^{Y}(\sigma_Y)\cdot Y.
\end{equation}
\item {\bf Even conormal index in codim 2}: If $codim(X)=2$ and $\sigma\in K^0(^bT^*X)$ we have that
\begin{equation}
\sum_{Y\in F_2(X)}Ind_{AS,2}^{Y}(\sigma_Y)\cdot Y\in Ker\,\delta_2=H_2^{cn}(X) 
\end{equation}
and 
\begin{equation}
Ind_{cn,ev}^X(\sigma)=\sum_{Y\in F_2(X)}Ind_{AS,2}^{Y}(\sigma_Y)\cdot Y.
\end{equation}
\item {\bf Odd conormal index in codim 2}: If $codim(X)=2$ and $\sigma\in K^1(^bT^*X)$ we have that
\begin{equation}
\sum_{Z\in F_1(X)}Ind_{Fred,1}^{\overline{Z}}(\sigma_{0,Z}^{nc})\cdot Z \in Ker\,\delta_1
\end{equation}
and 
\begin{equation}
Ind_{cn,odd}^X(\sigma)=\left[\sum_{Z\in F_1(X)}Ind_{Fred,1}^{\overline{Z}}(\sigma_{0,Z}^{nc})\cdot Z \right],
\end{equation}
where $[\cdot]$ denotes the class in $H_1^{cn}(X)$.
\item {\bf Even conormal index in codim 3}: If $codim(X)=3$ and $\sigma\in K^0(^bT^*X)$ we have that
\begin{equation}
\sum_{Z\in F_2^{cn}(X)}Ind_{Fred,2}^{\overline{Z}}(\sigma_{0,Z}^{nc})\cdot Z \in Ker\,\delta_2
\end{equation}
and 
\begin{equation}
Ind_{cn,ev}^X(\sigma)=\left[\sum_{Z\in F_2^{cn}(X)}Ind_{Fred,2}^{\overline{Z}}(\sigma_{0,Z}^{nc})\cdot Z \right],
\end{equation}
where $[\cdot]$ denotes the class in $H_2^{cn}(X)$.
\item {\bf Odd conormal index in codim 3}: If $codim(X)=3$ and $\sigma\in K^1(^bT^*X)$ we have that
\begin{equation}
\sum_{Y\in F_3^{cn}X)}Ind_{AS,3}^{Y}(\sigma_Y)\cdot Y \in Ker\,\delta_3,
\end{equation}

\begin{equation}
\sum_{Z\in F_1(X)}Ind_{Fred,1}^{\overline{Z}}(\sigma_{1,Z}^{nc})\cdot Z \in Ker\,\delta_1
\end{equation}
and 
\begin{equation}
Ind_{cn,odd}^X(\sigma)=(\sum_{Y\in F_3^{cn}(X)}Ind_{AS,3}^Y(\sigma_Y)\cdot Y,\left[\sum_{Z\in F_1(X)}Ind_{Fred,1}^{\overline{Z}}(\sigma_{1,Z}^{nc})\cdot Z \right])
\end{equation}
in $H_3^{odd}(X)=Ker\delta_3 \bigoplus H_1^{cn}(X)$,
where $[\cdot]$ denotes the class in $H_1^{cn}(X)$.
\end{itemize}
\end{theorem}

For the notations $F_\bullet^{cn}$ above and $F_\bullet^{cn,\delta}$ below see section \ref{sectionevenoddconormalindex}. These are both subsets of the set of faces of the given codimension, it is very interesting that not all the corners indices contribute in the computations above, we will come to this elsewhere.

Now, using the explicit computations of $Ind^X_{cn,ev}$ above we can specialize Theorem \ref{thmintroconormalobstruction} for codimensions 2 and 3 to get:

For codimension two:

\begin{corollary}[Fredholm anomalies in codimension 2]
Let $D$ be an elliptic $b$-pseudodifferential operator on a compact manifold with corners $X$ of codimension 2. Then  $D$ satisfies $(\cS\cF\cP)$ if and only if for every $Y\in F_2^{cn}(X)$ (face of codimension 2 
of $X$ belonging to a cycle, definition \ref{defF2cn}) the 2-suspended Atiyah-Singer index
\begin{equation}
Ind^Y_{AS,2}([\sigma_b(D)]_Y)=0.
\end{equation}
In particular, if $D$ satisfies $(\cF\cP)$ then all the 2-suspended Atiyah-Singer indices above vanish.
\end{corollary}

and, for codimension three:

\begin{corollary}[Fredholm anomalies in codimension 3]
Let $D$ be an elliptic $b$-pseudodifferential operator on a compact manifold with corners $X$ of codimension 3. Then  $D$ satisfies $(\cS\cF\cP)$ if and only if 
\begin{enumerate}
\item For every $Z\in F_2^{cn}(X)\setminus F_2^{cn,\delta}(X)$ the 2-suspended Fredholm index 
\begin{equation}
Ind^{\overline{Z}}_{Fred,2}([\sigma_b(D)]_{0,Z}^{nc})=0,
\end{equation}
and
\item there is a 3-conormal chain $a\in C_3^{cn}(X)$ such that
\begin{equation}
\sum_{Z\in F_2^{cn,\delta}(X)}Ind_{Fred,2}^{\overline{Y}}([\sigma_b(D)]_{0,Z}^{nc})\cdot Z=\delta_3(a).
\end{equation}
\end{enumerate}
In particular, if $D$ satisfies $(\cF\cP)$ then the following two points above hold.
\end{corollary}

{\bf Topological Fredholm anomalies:} For codimensions 2 and 3 the even conormal indices are computed in terms of suspended Atiyah-Singer and Fredholm morphisms respectively. In codimension 2, by the fundamental and by now classic Atiyah-Singer index theorem we can compute the Atiyah-Singer morphisms by a topological index formula and then the above even conormal index admits a computation in terms of these topological indices, see Corollary \ref{topevencncodim2}. In particular we have the following obstruction result:

\begin{theorem}[Topological anomalies in codimension 2]
Let $D$ be an elliptic $b$-pseudodifferential operator on a compact manifold with corners $X$ of codimension 2. Then  $D$ satisfies $(\cS\cF\cP)$ if and only if for every $Y\in F_2^{cn}(X)$ (face of codimension 2 
of $X$ belonging to a cycle, definition \ref{defF2cn}) the 2-suspended topological Atiyah-Singer index
\begin{equation}
\langle ch(\beta_2^{-1}([\sigma_b(D)]_Y))\wedge Td_Y, [Y]\rangle=0.
\end{equation}
In particular, if $D$ satisfies $(\cF\cP)$ then all the 2-suspended topological Atiyah-Singer indices above vanish.
\end{theorem}

For the codimension three case, the indices $Ind_{Fred,2}^{\overline{Z}}$ appearing in the computation of the even conormal index are the suspended version of the Fredholm indices for manifolds with boundary (mwc of codimension 1) and these indices were computed topologically by a cohomological formula in \cite{CLM}, as recalled in 
section \ref{subsectiontopanomalies3}. In particular, the above even conormal index admits a computation in terms of these topological indices, see corollary \ref{topevencncodim3}, and we have the following obstruction result where we are denoting by $Ind_{topFred,2}^{\overline{Y}}$ the suspended topological Fredholm indices defined in \cite{CLM} and recalled in \ref{deftopfredindex} above.

\begin{theorem}[Topological corner anomalies in codimension 3]
Let $X$ be a closed connected manifold with corners of codimension 3. Then an elliptic $b$-operator $P\in \overline{\Psi^0}(\Gamma_b(X))$ has the $(SFP)$ if and only if 
\begin{enumerate}
\item For every $Z\in F_2^{cn}(X)\setminus F_2^{cn,\delta}(X)$ the 2-suspended topological Fredholm index 
\begin{equation}
Ind^{\overline{Z}}_{topFred,2}([\sigma_b(D)]_{0,Z}^{nc})=0,
\end{equation}
and
\item there is a 3-conormal chain $a\in C_3^{cn}(X)$ such that
\begin{equation}
\sum_{Z\in F_2^{cn,\delta}(X)}Ind_{topFred,2}^{\overline{Z}}([\sigma_b(D)]_{0,Z}^{nc})\cdot Z=\delta_3(a). 
\end{equation}
\end{enumerate}
In particular, if $P$ satisfies $(\cF\cP)$ then the following two points above hold.
\end{theorem}

{\bf Final questions and remarks:}
\begin{enumerate}
\item The first thing we want to settle down is about the question of lower codimensional manifolds with corners. As mentioned above, in \cite{CarLes} we computed the isomorphisms (\ref{Tisointro}) for lower codimensions in terms of integral conormal homology using the computations of the first differentials of the spectral sequence canonically associated to the natural filtration by faces of a manifold with corners. As we dicussed in that paper the reason why we could not go further in the computation, with the methods of that paper, was the possibility that the conormal homology groups contained torsion but also the difficulty of computing higher differentials for the spectral sequence. In that paper we proved that the first conormal groups $H_1^{cn}(X)$ are always free and even if we did not know about the torsion from $H_2^{cn}$ we could perform the computation (\ref{Tisointro}) up to codimensional three manifold with corners. The question of torsion or not remained open. Later on, in \cite{CRLV} in collaboration with Mario Velasquez, we constructed the morphism (\ref{Trationalintro}) for any codimension and proved that it is rational isomorphism. The methods used in this last reference are completely different, indeed, in \cite{CRLV} we constructed a topological space whose topological K-theory is isomorphic to the $K-$theory of the $C^*$-algebra $\cK_b(X)$ by a Connes-Thom isomorphism, then we proved that the periodic singular cohmology of this space is canonically isomorphic to the periodic conormal homology groups. Now, the question of wether or not we could continue to use our integral coefficient methods was completely solved in \cite{SV} where Schick and Velasquez showed that from $H_2^{cn}$ there could be torsion but even further that the computations of the higher differentials in our first paper is as difficult as computing the higher differentials in topological K-theory for $CW$-complexes. In fact, for any manifold with (embedded) corners $X$ they associate to it a finite simplicial complex $\Sigma_X$ and show that its reduced $K$-homology is isomorphic to the $K-$theory group $K_*(\mathcal{K}_b(X))$ (up to a shift of degree). Also, one of their principal results is that for an arbitrary finite simplicial complex $\Sigma$ there is a manifold with corners $X$ such that $\Sigma\cong \Sigma_X$. As a consequence, the homology and $K-$homology which occur for finite simplicial complexes also occur as conormal homology groups of manifolds with corners. In particular, these groups can contain torsion. 

\item The last point justifies our choice of study and explicit the computations in integral coefficients homology for lower codimensions only. In a future work we will consider the even/odd rational conormal indices obtained by composing the analytic index morphism of a manifold with corners with the morphism (\ref{Trationalintro}). We will see that the classifying space used in \cite{CRLV} is part of a classifyng space giving the rational computation of the entire rational conormal morphisms.
\end{enumerate}
\section{Around manifolds with corners}

We start by defining the manifolds with corners we will be using in the entire paper.

A manifold with corners is a Hausdorff space covered by compatible coordinate charts with coordinate functions modeled in the spaces
$$\bR_k^n:=[0,+\infty)^k\times \bR^{n-k}$$
for fixed $n$ and possibly variable $k$.

\begin{definition}
A manifold with embedded corners $X$ is a Hausdorff topological space endowed with a subalgebra $C^{\infty}(X)\in C^0(X)$ satisfying the following conditions:
\begin{enumerate}
\item there is a smooth manifold $\tilde{X}$ and a map $\iota:X\to \tilde{X}$ such that $$\iota^*(C^\infty(\tilde{X}))=C^\infty(X),$$
\item there is a finite family of functions $\rho_i\in C^\infty(\tilde{X})$, called the defining functions of the hyperfaces, such that
$$\iota(X)=\bigcap_{i\in I}\{\rho_i\geq0\}.$$
\item for any $J\subset I$, 
\begin{center}
$d_x\rho_i(x)$ are linearly independent in $T^*_x\tilde{X}$ for all $x\in F_J:=\bigcap_{i\in J}\{\rho_i=0\}$.
\end{center}
\end{enumerate}
\end{definition}

{\bf Terminology:} In this paper we will only be considering manifolds with embedded corners. We will refer to them simply as manifolds with corners. We will also assume our manifolds to be connected. More general manifold with corners deserve attention but as we will see in further papers it will be more simple to consider them as stratified pseudomanifolds and desingularize them as manifolds with embedded corners with an iterated fibration structure.

The space  $X$ is naturally filtrated. Indeed, denote by $F_p$ the set of connected faces of codimension $p$ (and $d$ the codimension of $X$).  For a given face $f\in F_p$,  we define the index set $I(f)$ of $f$ to be the unique tuple  $(i_1,\ldots,i_p)$  such that $1\le i_1 <\ldots < i_p \le n$ and 
\begin{equation}
  f \subset H_{i_1}\cap\ldots\cap H_{i_p}
\end{equation}
where we recall that $H_j = \rho_j^{-1}(\{0\})\subset X$.  The filtration of $X$ is then given by:
\begin{equation}\label{filtration-by-codimension}
  X_{j} = \bigcup_{\stackrel{f\in F}{ \  d-j \le \mathrm{codim}(f)\le d}} f
\end{equation}
Indeed, we obtain the following filtration by closed subspaces:
\begin{equation}
 F_d = X_{0} \subset X_{1} \subset \cdots \subset X_d = X.
\end{equation}

\subsection{Conormal homology for manifolds with corners}\label{subsectionconormalhomology}
 Conormal homology is introduced  (under a different name) and studied in  \cite{Bunke}. In \cite{CarLes}, a slighty different presentation of this homology is given, after the observations that it  coincides with the $E^2$ page of the spectral sequence computing $K_*(C^*(\Gamma_b(X)))$ and that it should provide easily computable obstructions to various Fredholm perturbations properties. We just briefly recall the definition of the chain complex and of the differential of conormal homology (see \cite{Bunke,CarLes} for more details). 
 
 With the same notation as above, the  chain complex $C_*(X)$  is the $\ZZ$-module where $C_p(X)$  is generated by 
 \begin{equation}
 \{ f\otimes \varepsilon \ ;\ f\in F_p \text{ and } \varepsilon \text{ is an orientation of } N_f   \}.
\end{equation}
Here   $N_f = (T_fX/T f)^*$ is the conormal bundle of $f\subset X$. Note that this bundle is always trivializable with $e_i=d\rho_i, \ i\in I(f)$ as a preferred global basis, and oriented by $\epsilon_{f}:= \wedge_{i\in I(f)}e_i$ or its opposite.  
We define the differential  $\delta_* : C_*(X)\to C_{*-1}(X)$ by 
\begin{equation}\label{diffpcn1}
   \delta_p(f\otimes \epsilon) = \sum_{\substack{g\in F_{p-1},  \\ f\subset\overline{g}}} g\otimes e_{i_{(g,f)}}\lrcorner\epsilon
\end{equation}
where the index $i(f,g)$ corresponds to the defining function $\rho_{i(f,g)}$ defining $f$ in $g$ and where $\lrcorner$ denotes the contraction of exterior forms. The {\it conormal homology} of $X$, denoted by $H^{\cn}_*(X)$ is defined to be the homology of $(C_*(X),\delta_*)$. Even and odd groups are called  {\it periodic conormal homology}:
 \begin{equation}
 H^{\pcn}_{0}(X)=\oplus_{p\ge 0} H^{\cn}_{2p}(X) \text{ and }H^{\pcn}_{1}(X)=\oplus_{p\ge 0} H^{\cn}_{2p+1}(X).
\end{equation}
We can consider conormal homology with rational coefficients as well.
\subsection{K-theory of the $b-$groupoid of a mwc and conormal homology}\label{subsectionKvsconormal}

Let $X$ be a manifold with embedded corners, so by definition we are assuming there is a smooth manifold (of the same dimension) $\tilde{X} $ with $X\subset \tilde{X}$ and $\rho_1,..., \rho_n$ defining functions of the faces. In \cite{Mont}, Monthubert constructed a Lie groupoid (called Puff groupoid) associated to any decoupage $(\tilde{X}, (\rho_i))$, it has the following expression
\begin{equation}\label{Puffgrpd}
G(\tilde{X},(\rho_i))= \{(x,y,\lambda_1,...,\lambda_n)\in \tilde{X}\times \tilde{X}\times \RR^n: \rho_i(x)=e^{\lambda_i}\rho_i(y)\}.
\end{equation}
as a Lie subgroupoid of $\tilde{X}\times \tilde{X}\times \RR^k$. The Puff groupoid is not s-connected, denote by $G_c(\tilde{X},(\rho_i))$ its s-connected component.


\begin{definition}[The $b$-groupoid]\label{def:b-groupoid}
The $b-$groupoid $\Gamma_b(X)$ of $X$ is by definition the restriction to $X$ of the s-connected Puff groupoid (\ref{Puffgrpd}) considered above, that is
\begin{equation}
\Gamma_b(X):= G_c(\tilde{X},(\rho_i))|_X\rightrightarrows X
\end{equation}
\end{definition}

The $b-$groupoid was introduced by B. Monthubert in order to give a groupoid description for the Melrose's algebra of $b$-pseudodifferential operators, we will come to this in the last section. In any case, as a $C^{\infty,0}-$groupoid the $b$-groupoid is amenable, in particular its reduced and maximal $C^*$-algebras coincide, see \cite{Mont}. In this paper we will simply denote by $\cK_b(X)=C^*(\Gamma_b(X))$ this algebra and by
$$K^*(\Gamma_b(X)):=K_*(\cK_b(X))$$ 
its associated K-theory groups.

In this section we will recall the group isomorphisms
\begin{equation}\label{IsoCRL}
\xymatrix{
K^*(\Gamma_b(X))\ar[r]^-T_-\cong & H^{pcn}_*(X)
} 
\end{equation}
we gave in \cite{CarLes} for mwc $X$ of codimension up to 3 for $*=0,1$ (even, odd). The presentation we will give below slightly differs from the original one presented in \cite{CarLes}, in fact in ref.cit. the isomorphisms were given from the periodic conormal groups to the K-theory groups. We will give some detailed review since we will use it in the sequel.

{\bf Notation/Warning:} We use below the natural filtration of a given mwc $X$ introduced in \ref{filtration-by-codimension} above. But also, it will be convenient to use the associated filtration of algebras that has an indexation in the opposite direction, explicitly, we let, for $k\in \{0,...,d\}$ (where  $d=codim(X)$)
\begin{equation}
A_k:=C^*(\Gamma_b(X)|_{X_d\setminus X_{d-(k+1)}})
\end{equation}
with the convention $X_{-1}=\emptyset$. Then we have a filtration by $C^*$-ideals ($X_d\setminus X_{d-(k+1)}$ are open in $X$)
\begin{equation}
A_0\subset A_1\subset \cdots A_d=C^*(\Gamma_b(X))=\cK_b(X).
\end{equation}

{\bf 1A. $Codim(X)=0$:} The only face of codimension 0 is $X$ (we are always assuming $X$ to be connected). The isomorphism 
\begin{equation}
\xymatrix{
K_0(A_0)\ar[r]^-{T_0}_-\cong & H_0^{cn}(X)
}
\end{equation}
is simply given by sending the rank one projector $p_{X}$ chosen above to $X$.  

{\bf 1B. $Codim(X)=1$:} Consider the canonical short exact sequence
$$
\xymatrix{
0\ar[r]&A_0\ar[r]&A_1\ar[r]&A_1/A_0\ar[r]&0
}
$$
That gives the following exact sequence in K-theory
$$
\xymatrix{
0\ar[r]&K_1(A_1)\ar[r]&K_1(A_1/A_0)\ar[r]^-{d_1}&K_0(A_0)\ar[r]&0
}
$$
from which $K_1(A_1)\cong \ker d_1$. By theorem 5.27 and corollary 5.35 in \cite{CarLes}, we have the following commutative diagram
$$
\xymatrix{
K_1(A_1/A_0)\ar[d]_-{T_{1,0}}^-{\cong}\ar[r]^-{d_1}&K_0(A_0)\ar[d]^-{T_0}_-{\cong}\\
C_1^{cn}(X)\ar[r]_-{\delta_1}&C_0^{cn}(X).
}
$$

Then there is a unique isomorphism 
\begin{equation}
\xymatrix{
K_1(\Gamma_b(X))=K_1(A_1)\ar[r]^-{T_1}_-\cong & Ker\,\delta_1=H_1^{pcn}(X)
}
\end{equation}
fitting the following commutative diagram
$$
\xymatrix{
0\ar[r]&K_1(A_1)\ar[r]\ar[d]_-{T_1}^-\cong & K_1(A_1/A_0)\ar[d]_-{T_{1,0}}^-\cong\ar[r]^-{d_1}&K_0(A_0)\ar[d]_-{T_0}^-\cong\ar[r]&0\\
0\ar[r]&Ker\,\delta_1\ar[r]&C_1^{cn}(X)\ar[r]_-{\delta_1}&C_0^{cn}(X)\ar[r]&0.
}
$$

{\bf 1C. $Codim(X)=2$:} Let us start by defining the isomorphism for the even group $K_0(A_2)$. In fact, as shown in proposition 5.40 in \cite{CarLes} the morphism
\begin{equation}
\xymatrix{
K_0(A_2)\ar[r]^-{r}&K_0(A_2/A_0)
}
\end{equation}
induced by restriction to the boundary is an isomorphism. Now, let us consider the canonical short exact sequence
$$
\xymatrix{
0\ar[r]&A_1/A_0\ar[r]&A_2/A_0\ar[r]&A_1/A_0\ar[r]&0
}
$$
that gives the following exact sequence in K-theory
$$
\xymatrix{
0\ar[r]&K_0(A_2/A_0)\ar[r]&K_0(A_2/A_1)\ar[r]^-{d_2}&K_1(A_1/A_0)\ar[r]&0
}
$$
from which $K_0(A_2/A_0)\cong \ker d_2$. Again, by 5.27 a,d 5.35 in \cite{CarLes}, we have the following commutative diagram
$$
\xymatrix{
K_0(A_2/A_1)\ar[d]_-{T_{2,1}}^-{\cong}\ar[r]^-{d_2}&K_1(A_1/A_0)\ar[d]^-{T_{1,0}}_-{\cong}\\
C_2^{cn}(X)\ar[r]_-{\delta_2}&C_1^{cn}(X).
}
$$

Then there is a unique isomorphism 
\begin{equation}
\xymatrix{
K_0(\Gamma_b(X)|{X_1})=K_0(A_2/A_0)\ar[r]^-{T_2}_-\cong & Ker\,\delta_2=H_2^{cn}(X)
}
\end{equation}
fitting the following commutative diagram
$$
\xymatrix{
0\ar[r]&K_0(A_2/A_1)\ar[r]\ar[d]_-{T_2}^-\cong & K_0(A_2/A_1)\ar[d]_-{T_{2,1}}^-\cong\ar[r]^-{d_2}&K_1(A_1/A_0)\ar[d]_-{T_{1,0}}^-\cong\ar[r]&0\\
0\ar[r]&Ker\,\delta_2\ar[r]&C_2^{cn}(X)\ar[r]_-{\delta_2}&C_1^{cn}(X)\ar[r]&0.
}
$$
Let us now recall the isomorphism $T_2: K_1(A_2)\to H_1^{pcn}(X)$ for $X$ of codimension 2 in odd K-theory. In this case 
$$H_1^{pcn}(X)=H_1^{cn}(X)=\frac{Ker\,\delta_1}{Im\,\delta_2}.$$
Consider the following commutative diagram of short exact sequences
\begin{equation}
\xymatrix{
0\ar[r]&0\ar[r]&A_2/A_1\ar[r]&A_2/A_1\ar[r]&0\\
0\ar[r]&A_0\ar[u]\ar[r]&A_2\ar[u]\ar[r]&A_2/A_0\ar[u]\ar[r]&0\\
0\ar[r]&A_0\ar[u]\ar[r]&A_1\ar[u]\ar[r]&A_1/A_0\ar[u]\ar[r]&0.
}
\end{equation}
By applying K-theory to it we get (at least the part that concerns us here)
\begin{equation}
\xymatrix{
&K_0(A_2/A_1)\ar[d]\ar[r]^-{Id}_-=&K_0(A_2/A_1)\ar[d]^-{d_2}&\\
0\ar[r]&K_1(A_1)\ar[r]\ar[d]&K_1(A_1/A_0)\ar[d]\ar[r]^-{d_1}\ar[d]&K_0(A_0)\ar[d]^-{Id}_-=\\
0\ar[r]&K_1(A_2)\ar[d]\ar[r]&K_1(A_2/A_0)\ar[d]\ar[r]&K_0(A_0)\\
&0&0&
}
\end{equation}
As shown in \cite{CarLes} p.554-555 the isomorphism 
\begin{equation}
\xymatrix{
K_1(A_2)\ar[r]^-{T_2}& \frac{Ker\,\delta_1}{Im\,\delta_2}
}
\end{equation}
is then obtained, by a classic chase diagram argument using the previous diagram, as the unique isomorphism fitting the following commutative diagram
\begin{equation}
\xymatrix{
0&0\\
K_1(A_2)\ar[u]\ar[r]^-{T_2}_-\cong & \frac{Ker\,\delta_1}{Im\,\delta_2}\ar[u]\\
K_1(A_1)\ar[u]\ar[r]^-{T_1}_-\cong& Ker\,\delta_1\ar[u]_-q\\
K_0(A_2/A_1)\ar[u]\ar[r]_-{T_{2,1}}^-\cong&C_2^{cn}(X)\ar[u]_-{\delta_2}
}
\end{equation}
where $q:Ker\delta_1 \to \frac{Ker\delta_1}{Im\delta_2}$ denotes the canonical quotient projection.

{\bf 1D. $Codim(X)=3$:} Let us start by defining the isomorphism for the even group $K_0(A_3)$. Again, as shown in proposition 5.40 in \cite{CarLes} the morphism
\begin{equation}
\xymatrix{
K_0(A_3)\ar[r]^-{r}&K_0(A_3/A_0)
}
\end{equation}
induced by restriction to the boundary is an isomorphism. Consider the following commutative diagram of short exact sequences
\begin{equation}
\xymatrix{
&&0&0&\\
0\ar[r]&0\ar[r]&A_3/A_2\ar[u]\ar[r]&A_3/A_2\ar[r]\ar[u]&0\\
0\ar[r]&A_1/A_0\ar[u]\ar[r]&A_3/A_0\ar[u]\ar[r]&A_3/A_1\ar[u]\ar[r]&0\\
0\ar[r]&A_1/A_0\ar[u]\ar[r]& A_2/A_0 \ar[u]\ar[r]&A_2/A_1\ar[u]\ar[r]&0\\
&0\ar[u]&0\ar[u]&0\ar[u]&.
}
\end{equation}
By applying K-theory to it we get (at least the part that concerns us here)
\begin{equation}
\xymatrix{
&K_1(A_3/A_2)\ar[d]\ar[r]^-{Id}_-=&K_1(A_3/A_2)\ar[d]^-{d_3}&\\
0\ar[r]&K_0(A_2/A_0)\ar[r]\ar[d]&K_0(A_2/A_1)\ar[d]\ar[r]^-{d_2}\ar[d]&K_1(A_1/A_0)\ar[d]^-{Id}_-=\\
0\ar[r]&K_0(A_3/A_0)\ar[d]\ar[r]&K_0(A_3/A_1)\ar[d]\ar[r]&K_1(A_1/A_0)\\
&0&0&.
}
\end{equation}
As shown in \cite{CarLes} p.556-557 the isomorphism 
\begin{equation}
\xymatrix{
K_0(A_3/A_0)\ar[r]^-{T_3}& \frac{Ker\,\delta_2}{Im\,\delta_3}
}
\end{equation}
is then obtained, by a classic chase diagram argument using the previous diagram. In fact, to be more precise one gets, using the previous diagram, isomorphisms
$$K_0(A_2/A_0)\cong Ker\,\delta_2,$$
and
$$Im(K_1(A_3/A_2)\longrightarrow K_0(A_2/A_0))\cong Im\,\delta_3.$$
This will be important in the sequel.
Let us now recall the isomorphism $T_3: K_1(A_3)\to H_1^{pcn}(X)$ for $X$ of codimension 3 in odd K-theory. For this consider the following commutative diagram of short exact sequences
\begin{equation}
\xymatrix{
&0\ar[d]&0\ar[d]&&\\
0\ar[r]&A_1\ar[r]\ar[d]&A_1\ar[d]\ar[r]&0\ar[d]\ar[r]&0\\
0\ar[r]&A_2\ar[d]\ar[r]&A_3\ar[d]\ar[r]&A_3/A_2\ar[d]\ar[r]&0\\
0\ar[r]&A_2/A_1\ar[d]\ar[r]&A_3/A_1\ar[d]\ar[r]&A_3/A_2\ar[d]\ar[r]&0\\
&0&0&0&.
}
\end{equation}
By applying $K-$theory to the above diagram we get
\begin{equation}
\xymatrix{
0\ar[r]&K_1(A_1)\ar[r]^-=\ar[d]&K_1(A_1)\ar[d]\ar[r]&0&\\
0\ar[r]&K_1(A_2)\ar[d]\ar[r]&K_1(A_3)\ar[d]\ar[r]&K_1(A_3/A_2)\ar[d]^-= &\\
&0\ar[r]&K_1(A_3/A_1)\ar[r]&K_1(A_3/A_2)\ar[r]^-{d_3}&\\
}
\end{equation}
As shown in \cite{CarLes} p.557, up to a choice of a splitting of $K_1(A_3)\to K_1(A_3/A_1)$ (which is possible since $K_1(A_3/A_1)$ is free as a $\mathbb{Z}-$submodule of $K_1(A_3/A_2)$), there is an isomorphism (that depends of course on the choice of the splitting) between $K_1(A_3)$ and the direct sum
$$K_1(A_3/A_1)\bigoplus \frac{K_1(A_1)}{Im(K_0(A_2/A_1)\longrightarrow K_1(A_1))},$$
that corresponds to the isomorphism between $K_1(A_3)$ and the direct sum
$$Ker\,\delta_3 \bigoplus \frac{Ker\, \delta_1}{Im\,\delta_2}=H^{pcn}_{odd}(X).$$

\subsection{Index theory morphisms for manifolds with corners}

\subsubsection{Even/Odd analytic index morphism of a mwc}

Let $X$ be a (connected) mwc. We defined above its $b-$groupoid
$$\Gamma_b(X)\rightrightarrows X$$
as the restriction to $X$ of a Lie groupoid $G_c(\tilde{X},(\rho_i))\rightrightarrows \tilde{X}$. As for any Lie groupoid we can consider the associated tangent groupoid 
$$G_c^{tan}(\tilde{X},(\rho_i))\rightrightarrows \tilde{X}\times [0,1]$$
and its restriction to $X\times [0,1]$
\begin{equation}\label{tangentgrpd}
\Gamma_b^{tan}(X)\rightrightarrows X\times [0,1]
\end{equation}
which is what we call the tangent groupoid of the groupoid $\Gamma_b(X)$.
One of its main features is that its $C^*$-algebra fits in the following short exact sequence (of maximal $C^*-$algebras)
\begin{equation}
\xymatrix{
0\ar[r]&C^*(\Gamma_b(X)\times (0,1])\ar[r]&C^*(\Gamma_b^{tan}(X))\ar[r]^-{ev_0}&C^*(^bTX)\ar[r]&0
}
\end{equation}
where $ev_0$ is the morphism induced by the restriction/evaluation at $t=0$. In particular, since the $C^*$-algebra on the left is contractible, because of the presence of the interval $(0,1]$, one has that the morphism (in even and odd $K$-theory) induced by the evaluation at zero
\begin{equation}
\xymatrix{
K^*(\Gamma_b^{tan}(X))\ar[r]^-{e_0}_-\cong & K^*_{top}(^bT^*X)
}
\end{equation}
is an isomorphism, where at the right we can consider that we have the topological K-theory group of the space $^bT^*X$ since the convolution algebra $C^*(^bTX)$ and punctual product algebra $C_0(^bT^*X)$ are isomorphic by the fiberwise Fourier transform isomorphism.

For $X$ a mwc we define the even/odd analytic index morphism (with $*=0 \text{ or } 1$ if ev or odd respectively)
\begin{equation}
\xymatrix{
K^*(^bT^*X)\ar[rr]^-{Ind^X_{an,ev/odd}}&&K^*(\Gamma_b(X))
}
\end{equation}
as the morphism fitting (and hence defined by it) in the following commutative diagram
\begin{equation}
\xymatrix{
K^*(^bT^*X)\ar[rr]^-{Ind^X_{an,ev/odd}}&&K^*(\Gamma_b(X))\\
&K^*(\Gamma_b^{tan}(X))\ar[lu]^-{e_0}_-\cong \ar[ru]_-{e_1}.&
}
\end{equation}

In the case $X$ is closed we will recall in \ref{subsectionAnalyticindexmorphism} that these indices coincide with their analytically defined version, hence the name.

\subsubsection{The Fredholm index morphism of a mwc}\label{subsectionFredindexgeo}

For the analytic indices above we used the deformation tangent groupoid, for defining the so-called Fredholm index morphism we will need to recall the so-called Fredholm groupoid which is as well a deformation groupoid whose deformation index, as we will see below computes the Fredholm index.

Hence, for $X$ a mwc the Fredholm groupoid is defined as the open subgroupoid of the tangent groupoid
\begin{equation}\label{Fredgroupoid}
\Gamma_b^{Fred}(X):=\Gamma_b^{tan}(X)|_{X_F}\rightrightarrows X_F
\end{equation}
where $X_F:=(X\times [0,1])\setminus \partial X \times \{1\}$ (which is open in $X\times [0,1]$). The Fredholm groupoid contains 
$$\intx\times \intx \times (0,1]\rightrightarrows \intx\times (0,1]$$ as an open saturated subgroupoid whose closed complement groupoid writes down as
\begin{equation}
^bT_{nc}X:=\Gamma_b^{Fred}(X)|_{X_\partial}\rightrightarrows X_\partial
\end{equation}
where $X_\partial:=X_F\setminus (\intx\times (0,1])$ and it is called {\it the noncommutative tangent space of $X$}. In fact there is a short exact sequence
\begin{equation}
\xymatrix{
0\ar[r]&C^*(\intx\times \intx \times (0,1])\ar[r]&C^*(\Gamma_b^{Fred}(X))\ar[r]^-{r_{X_\partial}}&C^*(^bT_{nc}X)\ar[r]&0
}
\end{equation}
where $r_{X_\partial}$ is the morphism induced by the restriction to $X_\partial$. In particular, since the $C^*$-algebra on the left is contractible, because of the presence of the interval $(0,1]$, one has that the morphism (in even and odd $K$-theory) induced by the restriction to $X_\partial$
\begin{equation}
\xymatrix{
K^*(\Gamma_b^{Fred}(X))\ar[r]^-{r_\partial}_-\cong & K^*(^bT_{nc}X)
}
\end{equation}
is an isomorphism.

For $X$ a mwc we define the Fredholm index morphism 
\begin{equation}
\xymatrix{
K^0(^bT_{nc}X)\ar[rr]^-{Ind^X_{Fred}}&&K^0(\intx \times \intx)\cong \mathbb{Z}
}
\end{equation}
as the morphism fitting (and hence defined by it) in the following commutative diagram
\begin{equation}
\xymatrix{
K^0(^bT_{nc}X)\ar[rr]^-{Ind^X_{Fred}}&&K^0(\intx \times \intx)\cong \mathbb{Z}\\
&K^0(\Gamma_b^{Fred}(X))\ar[lu]^-{r_\partial}_-\cong \ar[ru]_-{e_1}.&
}
\end{equation}

For $X$ closed we will recall in \ref{Fredsubsection} that this index coincide with its analytically defined version that computes the Fredholm index for fully elliptic operators. In \cite{CLM}, for the case of manifolds with boundaries, this index morphism was called the APS index morphism associated to such a manifold, the name was inspired by the foundational work of Atiyah, Patodi and Singer \cite{APS1} and its sequels.

\subsubsection{Even/Odd conormal index morphism of a mwc}

For a mwc $X$ we define the even/odd conormal index morphism (with $*=0 \text{ or } 1$ if ev or odd respectively)
\begin{equation}
\xymatrix{
K^*(^bT^*X)\ar[rr]^-{Ind^X_{cn,ev/odd}}&&H_*^{pcn}(X)
}
\end{equation}
given as the composition of the even/odd analytic index morphism followed by the isomorphism
\begin{equation}
\xymatrix{
K^*(\Gamma_b(X))\ar[r]^-T_-\cong & H_*^{pcn}(X)
}
\end{equation}

\subsubsection{Suspended Atiyah-Singer index morphisms}

Let $Y$ be a smooth manifold. For $p\in \mathbb{N}$ we define the suspended (or $p$-suspended if needed) Atiyah-Singer index morphism associated to $Y$ to be the morphism ($*=0$ if $p$ even or $*=1$ if $p$ odd)
\begin{equation}
\xymatrix{
K^*(T^*Y\times \mathbb{R}_p)\ar[rr]^-{Ind_{AS,p}^Y}&& \mathbb{Z}
}
\end{equation}
where $\mathbb{R}_p$ stands above for the additive group $(\mathbb{R}^p,+)$, and given as the composition of the deformation index morphism $(e_1\circ e_0^{-1})$ given by the tangent groupoid $G_Y^{tan}\rightrightarrows Y\times [0,1]$, 
\begin{equation}
\xymatrix{
K^*(T^*Y\times \mathbb{R}_p)&K^*(G_Y^{tan}\times \mathbb{R}_p)\ar[r]^-{e_1}\ar[l]_-{e_0}^-\cong& K^*(Y\times Y\times \mathbb{R}_p),
}
\end{equation}
followed by the Bott inverse isomorphism
\begin{equation}
\xymatrix{
K^*(Y\times Y\times \mathbb{R}_p)\ar[r]^-{\beta_p^{-1}}&\mathbb{Z}
}
\end{equation}
with the Bott element $\beta_p$ considered in the introduction.
\subsubsection{Suspended Fredholm index morphisms}

Let $Z$ be a manifold with corners. For $p\in \mathbb{N}$ we define the suspended (or $p$-suspended if needed) Fredholm index morphism associated to $Z$ to be the morphism ($*=0$ if $p$ even or $*=1$ if $p$ odd)
\begin{equation}
\xymatrix{
K^*(^bT_{nc}Z\times \mathbb{R}_p)\ar[rr]^-{Ind_{Fred,p}^Z}&& \mathbb{Z}
}
\end{equation}
where $\mathbb{R}_p$ stands above for the additive group $(\mathbb{R}^p,+)$, and given as the composition of the deformation index morphism $(e_1\circ r_\partial^{-1})$ given by the Fredholm groupoid $\Gamma_b^{Fred}(Z)\rightrightarrows Z_F$, 
\begin{equation}
\xymatrix{
K^*(^bT_{nc}Z\times \mathbb{R}_p)&K^*(\Gamma_b^{Fred}(Z)\times \mathbb{R}_p)\ar[r]^-{e_1}\ar[l]_-{r_\partial}^-\cong& K^*(\overset{\:\circ}{Z}\times \overset{\:\circ}{Z}\times \mathbb{R}_p),
}
\end{equation}
followed by the Bott inverse isomorphism
\begin{equation}
\xymatrix{
K^*(\overset{\:\circ}{Z}\times \overset{\:\circ}{Z}\times \mathbb{R}_p)\ar[r]^-{\beta_p^{-1}}&\mathbb{Z}
}
\end{equation}
with the Bott element $\beta_p$ considered in the introduction.

\section{Computation of the Even and odd conormal index morphisms}\label{sectionevenoddconormalindex}

\subsection{Codimension 1}

Let $X$ be a mwc of codimension 1, we consider the odd conormal index morphism
\begin{equation}
\xymatrix{
K^1_{top}(^bT^*X)\ar[rr]^-{Ind_{cn,odd}^X}&&H_1^{pcn}(X)=Z_1^{cn}(X)
}
\end{equation}
given as the composition of the odd analytic index 
\begin{equation}
\xymatrix{
K^1_{top}(^bT^*X)\ar[rr]^-{Ind_{ana,odd}^X}&&K_1(C^*(\Gamma_b(X)))
}
\end{equation}
and the isomorphism
\begin{equation}
K_1(C^*(\Gamma_b(X))\stackrel{T}{\longrightarrow}Z_1^{cn}(X)
\end{equation}
recalled above. We will compute the above mentioned odd conormal morphism in the case $\# F_1(X)\geq 2$ otherwise the group $Z_1^{cn}(X)$ is zero.

Before stating a precise formula for the index morphism above we need some terminology. Let $\sigma\in K^1_{top}(^bT^*X)$, for $Y\in F_1(X)$ there is a class given canonically by restriction
\begin{equation}
\sigma_Y\in K^1(T^*Y\times \mathbb{R}_1),
\end{equation}
for which we can take its $1$-suspended Atiyah-Singer index $Ind_{AS,1}^{Y}(\sigma_Y)\in \mathbb{Z}$. We have the following theorem:


\begin{theorem}
Let $X$ be a closed connected manifold with boundary with two or more boundary components. Let  $\sigma\in K^1(^bT^*X)$, with the notations above we have that
\begin{equation}
\sum_{Y\in F_1(X)}Ind_{AS,1}^{Y}(\sigma_Y)\cdot (Y\otimes \epsilon_Y) \in Ker\,\delta_1
\end{equation}
and 
\begin{equation}
Ind_{cn,odd}^X(\sigma)=\sum_{Y\in F_1(X)}Ind_{AS,1}^{Y}(\sigma_Y)\cdot (Y\otimes \epsilon_Y).
\end{equation}
\end{theorem}

\begin{proof}
Let us briefly recall the computation of $K^1(\Gamma_b)=K^1(\Gamma_b|_{X_1})$ in terms of the conormal homology of the filtration $X_0\subset X_1$, indeed there are isomorphisms $T$ fitting in the following commutative diagram 
\begin{equation}
\xymatrix{
0\ar[r]&K^1(\Gamma_b|_{X_1})\ar[d]_-T^-\cong\ar[r]&K^1(\Gamma_b|_{X_0})\ar[d]_-T^-\cong\ar[r]&K^0(\Gamma_b|_{X_1\setminus  X_0})\ar[d]_-T^-\cong\\
0\ar[r]&Ker\,\delta_1 \ar[r]&C_1^{cn}(X)\ar[r]_-{\delta_1}&C_0^{cn}(X)
}
\end{equation}
where we can besides remark that the morphism $K^1(\Gamma_b|_{X_1})\longrightarrow K^1(\Gamma_b|_{X_0})$ corresponds exactly to the morphism 
$\oplus_{c\in F_1}r_c$ induced by the restriction to the faces of codimension 1.

The conclusion of the theorem follows from the following commutative diagram
\begin{equation}
\xymatrix{
K^1(\Gamma_b^{tan}(X))\ar[r]^-{e_1}\ar[d]_-{\oplus_c r_c}&K^1(\Gamma_b)\ar[d]^-{\oplus_c r_c}\ar[r]^-T_-\cong & Ker\,\delta_1 \ar[d]\\
\bigoplus_{Y\in F_1}K^1(G_Y^{tan}\times \mathbb{R}_1)\ar[r]_-{e_1}&\bigoplus_{Y\in F_1}K^1(Y\times Y\times \mathbb{R}_1)\cong K^0(\Gamma_b|_{X_0})\ar[r]_-T^-\cong&C_1^{cn}(X) 
}
\end{equation}
where the right vertical arrow corresponds to the inclusion $Ker\,\delta_1 \to C_1^{cn}(X)$.
\end{proof}

\subsection{Codimension 2}

\subsubsection{Even conormal index in codim 2}
Let $X$ be a mwc of codimension 2, we consider the even conormal index morphism
\begin{equation}
\xymatrix{
K^0_{top}(^bT^*X)\ar[rr]^-{Ind_{cn,ev}^X}&&H_0^{pcn}(X)=Z_2^{cn}(X)
}
\end{equation}
given as the composition of the boundary analytic index 
\begin{equation}
\xymatrix{
K^0_{top}(^bT^*X)\ar[rr]^-{Ind_\partial^X}&&K_0(C^*(\Gamma_b(X)|_{\partial(X)}))
}
\end{equation}
and the isomorphism
\begin{equation}
K_0(C^*(\Gamma_b(X)|_{\partial(X)}))\stackrel{T}{\longrightarrow}Z_2^{cn}(X)
\end{equation}
recalled above. Before stating a precise formula for the index morphism above we need some terminology. Let $\sigma\in K^0_{top}(^bT^*X)$, for $Y\in F_2(X)$ there is a class given canonically by restriction
\begin{equation}
\sigma_Y\in K^0(T^*Y\times \mathbb{R}_2)
\end{equation}
and an associated $2$-suspended Atiyah-Singer index $Ind_{AS,2}^{Y}(\sigma_Y)\in \mathbb{Z}$.

We need in this codimension 2 case an extra definition.

\begin{definition}\label{defF2cn}
Let $Y\in F_2(X)$. Consider the canonical coefficient projection
$$C_2^{cn}(X)\stackrel{q_Y}{\longrightarrow}\mathbb{Z}$$
that computes the coefficient of $Y\otimes \epsilon_Y$. We say that $Y$ belongs to a corner cycle if there is $m\in Z_2^{cn}(X)$ such that $q_Y(m)\neq 0$.
We denote by $F_2^{cn}(X)\subset F_2(X)$ the subset of codimension 2 faces belonging to a cycle.
\end{definition}

We have the following computation of the even conormal index morphism in terms of corner cycles.

\begin{theorem}\label{cnindexcodim2}
Let $[\sigma]\in K^0(^bT^*X)$, with the notations above, we have that
\begin{equation}
\sum_{Y\in F_2^{cn}(X)}Ind_{AS,2}^Y(\sigma_Y)\cdot (Y\otimes \epsilon_Y)\in Ker\,\delta_2
\end{equation}
and
\begin{equation}
Ind_{cn,ev}^X([\sigma])=\sum_{Y\in F_2^{cn}(X)}Ind_{AS,2}^Y(\sigma_Y)\cdot (Y\otimes \epsilon_Y).
\end{equation}
\end{theorem}

\begin{proof}
Let us briefly recall the computation of $K^0(\Gamma_b|_{\partial X})=K^0(\Gamma_b|_{X_1})=K_0(A_2/A_0)$ in terms of the conormal homology of the filtration $X_0\subset X_1\subset X_2$, indeed there are isomorphisms $T$ fitting in the following commutative diagram 
\begin{equation}
\xymatrix{
0\ar[r]&K^0(\Gamma_b|_{X_1})\ar[d]_-T^-\cong\ar[r]&K^0(\Gamma_b|_{X_0})\ar[d]_-T^-\cong\ar[r]&K^1(\Gamma_b|_{X_1\setminus  X_0})\ar[d]_-T^-\cong\\
0\ar[r]&Ker\,\delta_2 \ar[r]&C_2^{cn}(X)\ar[r]_-{\delta_1}&C_1^{cn}(X)
}
\end{equation}
where we can besides remark that the morphism $K^0(\Gamma_b|_{X_1})\longrightarrow K^0(\Gamma_b|_{X_0})$  corresponds exactly to the morphism 
$\oplus_{Y\in F_2}r_Y$ induced by the restriction to the faces of codimension 2.

The conclusion of the theorem follows from the following commutative diagram
\begin{equation}
\xymatrix{
K^0(\Gamma_b^{tan}(X))\ar[r]^-{e_1^\partial}\ar[d]_-{\oplus_Y r_Y}&K^0(\Gamma_b|_{\partial X})\ar[d]^-{\oplus r_Y}\ar[r]^-T_-\cong & Ker\,\delta_2 \ar[d]\\
\bigoplus_{Y\in F_2}K^0(G_Y^{tan}\times \mathbb{R}_2)\ar[r]_-{e_1}&\bigoplus_{Y\in F_2}K^0(Y\times Y\times \mathbb{R}_2)\cong K^0(\Gamma_b|_{X_0})\ar[r]_-T^-\cong&C_2^{cn}(X) 
}
\end{equation}
where the right vertical arrow corresponds to the inclusion $Ker\,\delta_2 \to C_2^{cn}(X)$.
\end{proof}

\subsubsection{Odd conormal index in codim 2}

Let $X$ be a mwc of codimension 2, we consider the odd conormal index morphism
\begin{equation}
\xymatrix{
K^1_{top}(^bT^*X)\ar[rr]^-{Ind_{cn,odd}^X}&&H_1^{pcn}(X)
}
\end{equation}
given as the composition of the odd analytic index 
\begin{equation}
\xymatrix{
K^1_{top}(^bT^*X)\ar[rr]^-{Ind_{ana,odd}^X}&&K_1(C^*(\Gamma_b(X)))
}
\end{equation}
and the isomorphism
\begin{equation}
K_1(C^*(\Gamma_b(X))\stackrel{T}{\longrightarrow} H_1^{cn}(X)=\frac{Ker\,\delta_1}{Im\,\delta_2}
\end{equation}
recalled above. Before stating a precise formula for the index morphism above we need some terminology. Let $Y\in F_1(X)$, in this case $\overline{Y}$ is a manifold with boundary. Consider the following commutative diagram
\begin{equation}
\xymatrix{
0&0\\
&\\
K^1(\Gamma_b^{tan}(X))\ar[r]^-{r_{\overline{Y}}}\ar[uu]^-{r_{X_0\times \{1\}}}&K^1(\Gamma_b^{tan}(\overline{Y})\times\mathbb{R}_1)\ar[uu]_-{r_{\partial \overline{Y}\times\{1\}}}\\
K^1(\Gamma_b^{tan}(X)|_{(X_0\times \{1\})^c})\ar[u]^-t\ar[r]_-{r_{\overline{Y}}}&K^1(\Gamma_b^{Fred}(\overline{Y})\times \mathbb{R}_1)\ar[u]_-{\iota_0}
}
\end{equation}
where the horizontal morphisms are induced from the indicated restrictions to closed subgroupoids, the columns are exact and induced from the indicated restrictions on the top, and where the first line is composed by zeros since $K^1(\Gamma_b^{tan}(X)|_{X_0\times \{1\}})=0$ and 
$K^1((\Gamma_b^{tan}(\overline{Y})\times\mathbb{R}_1)|_{\partial \overline{Y}\times\{1\}})=0 	$ in this case since $\overline{Y}$ is a closed manifold with boundary.

Let
\begin{equation}\label{fixsigma}
\sigma^{nc}_0\in K^1(\Gamma_b^{tan}(X)|_{(X_0\times \{1\})^c})
\end{equation}
be a fixed lifting of $e_0^{-1}(\sigma)\in K^1(\Gamma_b^{tan}(X))$ by the morphism $t$ above, {\it i.e.} $e_0(t(\sigma^{nc}_0))=\sigma\in K^1(^bT^*X)$. We let $\sigma^{nc}_{0,Y}\in K^1(T_{nc}\overline{Y} \times \mathbb{R}^1)$ to be the image of $(r_{\overline{Y}})(\sigma^{nc}_0) \in K^1(\Gamma_b^{Fred}(\overline{Y})\times \mathbb{R})$ under the isomorphism in K-theory
\begin{equation}
\xymatrix{
K^1(T_{nc}\overline{Y} \times \mathbb{R}_1)&K^1(\Gamma_b(\overline{Y})^{Fred}\times \mathbb{R}_1)\ar[l]^-{r_\partial}_-{\cong}.}
\end{equation}

We can hence consider the $1$-suspended Fredholm index $Ind_{Fred,1}^{\overline{Y}}(\sigma^{nc}_{0,Y})\in \mathbb{Z}$. We can now compute the odd conormal index in terms of these indices.

\begin{theorem}\label{Thmcnboundaryindex}
Let $\sigma\in K^1(^bT^*X)$, with the notations above we have that
\begin{equation}
\sum_{Y\in F_1(X)}Ind_{Fred,1}^{\overline{Y}}(\sigma^{nc}_{0,Y})\cdot (Y\otimes \epsilon_Y) \in Ker\, \delta_1
\end{equation}
and
\begin{equation}
Ind_{cn,odd}^X(\sigma)=\left[\sum_{Y\in F_1(X)}Ind_{Fred,1}^{\overline{Y}}(\sigma^{nc}_{0,Y})\cdot (Y\otimes \epsilon_Y)\right]
\end{equation}
where $[\cdot]$ denotes the class in $H_1^{cn}(X)$.
\end{theorem}
\begin{proof}
Let $\sigma^{tan}\in K^1(\Gamma_b^{tan}(X))$ with $e_0(\sigma^{tan})=\sigma\in K^1(^bT^*X)$. Denote 
$$J=e_1(\sigma^{tan})\in K_1(A_2).$$
We want to compute $J$ in $K_1(A_2)=K_1(A_1)/Im(K_0(A_2/A_1)\to K_1(A_1)) $  that we identified above with $H_1^{cn}(X)=Ker\,\delta_1/Im\, \delta_2$.
Consider the following commutative diagram with exact columns
\begin{equation}
\xymatrix{
0&0&0\\
&&\\
K^1(\Gamma_b^{tan}(X))\ar[r]^-{e_1}\ar[uu]^-{r_{X_0\times \{1\}}}&K_1(A_2)\ar[r]^-j\ar[uu]^-{r_{X_0}}&K_1(A_2/A_0)\ar[uu]^-{r_{X_0}}\\
K^1(\Gamma_b^{tan}(X)|_{(X_0\times \{1\})^c})\ar[u]^-t\ar[r]_-{e_1}&K_1(A_1)\ar[r]_-j\ar[u]^-{q}& K_1(A_1/A_0)\ar[u]^-{q}\\
&&
}
\end{equation}
By commutativity of the diagram above
$$q(I)=J\in K_1(A_2),$$
if we denote by $I=e_1(\sigma_{0}^{nc})\in K_1(A_1)=K^1(\Gamma_b(X)|_{X_2\setminus X_0})$.

Notice also that $I\in Ker(K_1(A_1/A_0)\stackrel{\delta_1}{\longrightarrow}K_0(A_0))$ and that
the projection $q:K_1(A_1)\longrightarrow K_1(A_2)$ corresponds in conormal homology to the canonical projection $Ker\delta_1 \to \frac{Ker\delta_1}{Im\delta_2}$. We would then prove the theorem if we can prove that
\begin{equation}
j(I)=(Ind_{Fred,1}^Y(\sigma^{nc}_{0,Y}))_{Y\in F_1}\in K_1(A_1/A_0)\cong \bigoplus_{Y\in F_1}K^1(Y\times Y \times \mathbb{R}_1),
\end{equation}
up to Bott isomorphism.

Consider the following commutative diagram
\begin{equation}
\xymatrix{
K^1(\Gamma_b^{tan}(X)|_{(X_0\times \{1\})^c})\ar[d]_-a\ar[r]^-{e_1}&K_1(A_1)\ar[r]^-j\ar[d]^-j&K_1(A_1/A_0)\\
K^1(\Gamma_b^{tan}(X)|_{(X_1\times [0,1])\setminus X_0\times \{1\}})\ar[r]_-F&K_1(A_1/A_0)\ar[ru]_-=&
}
\end{equation}
where $F:K^1(\Gamma_b^{tan}(X)|_{(X_1\times [0,1])\setminus X_0\times \{1\}})\to K^1(\Gamma_b(X)|_{X_1\setminus X_0})=K_1(A_1/A_0)$
is the morphism induced by evaluation at $\epsilon=1$.
In particular we have that
\begin{equation}
j(I)=F(a(\sigma_{0}^{nc}))\in K_1(A_1/A_0)\approx \bigoplus_{Y\in F_1}\mathbb{Z}\cdot Y
\end{equation}
We want then to show that the $Y-$component of $F(a(\sigma_0^{nc}))$, that we denote by $F(a(\sigma_0^{nc}))_Y$, is $Ind_{Fred,1}^Y(\sigma^{nc}_{0,Y})$, up to Bott isomorphism, for every $Y\in F_1(X)$.

Let us fix then $Y\in F_1(X)$, we have the following commutative diagram
\begin{equation}
\xymatrix{
K^1(\Gamma_b^{tan}(X)|_{(X_0\times \{1\})^c})\ar[r]^-a&K^1(\Gamma_b^{tan}(X)|_{(X_1\times [0,1])\setminus X_0\times \{1\}})\ar[d]_-{F}\ar[r]^-{r_{\overline{Y}}}&K^1(\Gamma_b^{tan}(X)|_{(\overline{Y}\times [0,1])\setminus \partial\overline{Y}\times \{1\}})\ar[d]^-{F_Y} \\
&K^1(\Gamma_b(X)|_{X_1\setminus X_0})\ar[r]_-{r_Y}&K^1(Y\times Y \times \mathbb{R}_1)
}
\end{equation}
where $r_{\overline{Y}}$ and $r_Y$ stand for the morphisms induced by the obvious restrictions indicated by the notations and where $F_Y$ is the morphism induced by evaluation at $\epsilon=1$.
We have just to remark that 
$$K^1(\Gamma_b^{tan}(X)|_{(\overline{Y}\times [0,1])\setminus \partial\overline{Y}\times \{1\}})=K^1(\Gamma_b^{Fred}(\overline{Y})\times \mathbb{R}_1)$$
to obtain
\begin{equation}\label{Fredc}
F(a(\sigma_0^{nc}))_Y=Ind_{Fred,1}(\sigma^{nc}_{0,Y}) 
\end{equation}
and get the theorem statement.
\end{proof}
\subsection{Codimension 3}

\subsubsection{Even conormal index in codim 3}
Let $X$ be a mwc of codimension 3, we consider the morphism
\begin{equation}
\xymatrix{
K^0_{top}(^bT^*X)\ar[rr]^-{Ind_\partial^X}&&H_0^{pcn}(X)=H_2^{cn}(X)
}
\end{equation}
given as the composition of the Boundary analytic index 
\begin{equation}
\xymatrix{
K^0_{top}(^bT^*X)\ar[rr]^-{Ind_\partial}&&K_0(C^*(\Gamma_b(X)|_{\partial(X)}))
}
\end{equation}
and the isomorphism (integral for this codimension)
\begin{equation}
K_0(C^*(\Gamma_b(X)|_{\partial(X)}))\stackrel{T}{\longrightarrow}H_2^{cn}(X).
\end{equation}
In the case the manifold $X$ is closed and connected we will give an analytic interpretation of the morphism above.

Let $\sigma\in K^{0}_{top}(^bT^*X)$, let $Y\in F_2(X)$ and consider $\overline{Y}$ which is a manifold with boundary. By restriction there is an induced symbol class 
$$\sigma_Y\in K^0(^bT^*\overline{Y}\times \mathbb{R}_2).$$
Now, by the isomorphism in $K-$theory
\begin{equation}
\xymatrix{
K^0(^bT^*\overline{Y}\times \mathbb{R}_2) & K^0(\Gamma_b(\overline{Y})^{tan}\times \mathbb{R}_2)\ar[l]^-{e_0}_-{\cong}
}
\end{equation}
and because
\begin{equation}
\xymatrix{
K_0(\Gamma_b(\overline{Y})^{Fred}\times \mathbb{R}_2)\ar[r]^-{i_0}&K^0(\Gamma_b(\overline{Y})^{tan}\times \mathbb{R}_2)\ar[r]&0
}
\end{equation}
is a surjection (because $\overline{Y}$ is a manifold with boundary), there are liftings of $\sigma_Y$ to $K_0(T_{nc}\overline{Y} \times \mathbb{R}_2)$ for the induced morphism $i_0$ fitting in the commutative diagram. 
\begin{equation}
\xymatrix{
K_0(T_{nc}\overline{Y} \times \mathbb{R}_2)\ar[r]^-{i_0}&K^0(^bT^*\overline{Y}\times \mathbb{R}_2)\\
K_0(\Gamma_b(\overline{Y})^{Fred}\times \mathbb{R}_2)\ar[r]_-{i_0}\ar[u]^-{r_\sigma}_-\cong &
K^0(\Gamma_b(\overline{Y})^{tan}\times \mathbb{R}_2)\ar[u]_-{e_0}^-\cong
}
\end{equation}
Now, not every lifting is useful for our purposes, we will need during the statement below particular liftings. Let us introduce these noncommutative symbols classes, for this consider the following commutative diagram
\begin{equation}
\xymatrix{
0&0&0\\
&&\\
K^0(\Gamma_b^{tan}(X))\ar[r]^-{r_1}\ar[uu]^-{r_{X_0\times \{1\}}}&K^0(\Gamma_b^{tan}(X)|_{X_1\times [0,1]})\ar[uu]^-{r_{X_0\times \{1\}}}\ar[r]^-{r_{\overline{Y}}}&K^0(\Gamma_b^{tan}(\overline{Y})\times\mathbb{R}_2)\ar[uu]^-{r_{\partial \overline{Y}\times\{1\}}}\\
K^0(\Gamma_b^{tan}(X)|_{(X_0\times \{1\})^c})\ar[u]^-t\ar[r]_-{r_1}&K^0(\Gamma_b^{tan}(X)|_{(X_1\times [0,1])\setminus (X_0\times\{1\})})
\ar[u]^-{\iota_0}\ar[r]_-{r_{\overline{Y}}}&K^0(\Gamma_b^{Fred}(\overline{Y})\times \mathbb{R}_2)\ar[u]^-{\iota_0}
}
\end{equation}
where the horizontal morphisms are induced from the indicated restrictions to closed subgroupoids, the columns are exact induced from the indicated restrictions on the top and where the first line is composed by zeros since $K^0(\Gamma_b^{tan}(X)|_{X_0\times \{1\}})=0$ and 
$K^0((\Gamma_b^{tan}(\overline{Y})\times\mathbb{R}_2)|_{\partial \overline{Y}\times\{1\}})=0 	$ in this case since $\overline{Y}$ is a closed manifold with boundary.

Let
\begin{equation}\label{fixsigma}
\sigma^{nc}_0\in K^0(\Gamma_b^{tan}(X)|_{(X_0\times \{1\})^c})
\end{equation}
be a fixed lifting of $e_0^{-1}(\sigma)\in K^0(\Gamma_b^{tan}(X))$ by the morphism $t$ above, {\it i.e.} $e_0(t(\sigma^{nc}_0))=\sigma\in K^0(^bT^*X)$. 
We let $\sigma^{nc}_{0,Y}\in K_0(T_{nc}\overline{Y} \times \mathbb{R}_2)$ to be the image of $(r_{\overline{Y}}\circ r_1)(\sigma^{nc}_0) \in K^0(\Gamma_b^{Fred}(\overline{Y})\times \mathbb{R}_2)$ under the isomorphism in K-theory
\begin{equation}
\xymatrix{
K_0(T_{nc}\overline{Y} \times \mathbb{R}_2)&K_0(\Gamma_b(\overline{Y})^{Fred}\times \mathbb{R}_2)\ar[l]^-{r_\partial}_-{\cong}.}
\end{equation}

We can consider the associated $2$-suspended Fredholm index  
$$Ind_{Fred,2}^{\overline{Y}}(\sigma^{nc}_{0,Y})\in \mathbb{Z}.$$


We can now compute the boundary analytic morphism in conormal homology in terms of these indices.

\begin{theorem}\label{Thmcnboundaryindex}
Let $\sigma\in K^0(^bT^*X)$, we have that
\begin{equation}
\sum_{Y\in F_2^{cn}(X)}Ind_{Fred,2}^{\overline{Y}}(\sigma^{nc}_{0,Y})\cdot (Y\otimes \epsilon_Y) \in Ker\, \delta_2
\end{equation}
and
\begin{equation}
Ind_{cn}^X(\sigma)=\left[\sum_{Y\in F_2^{cn}(X)}Ind_{Fred,2}^{\overline{Y}}(\sigma^{nc}_{0,Y})\cdot (Y\otimes \epsilon_Y)\right]
\end{equation}
where $[\cdot]$ denotes the class in $H_2^{cn}(X)$.
\end{theorem}
\begin{proof}
Let $\sigma^{tan}\in K^0(\Gamma_b^{tan}(X))$ with $e_0(\sigma^{tan})=[\sigma_P]\in K^0(^bT^*X)$. Denote 
$$J=e_1^\partial(\sigma^{tan})\in K_0(A_3/A_0),$$
where $e_1^\partial$ denotes the morphism induced by the evaluation at $\epsilon=1$ followed by the morphism induced by restriction to the boundary (which is an isomorphism in this case).

We want to compute $J$ in $K_0(A_3/A_0)=K_0(A_2/A_0)/Im(K_1(A_3/A_2)\to K_0(A_2/A_0)) $  that we identified above with $H_2^{cn}(X)=Ker\,\delta_2/Im\, \delta_3$.
Consider the following commutative diagram with exact columns
\begin{equation}
\xymatrix{
0&0&0\\
&&\\
K^0(\Gamma_b^{tan}(X))\ar[r]^-{e_1^\partial}\ar[uu]^-{r_{X_0\times \{1\}}}&K_0(A_3/A_0)\ar[r]^-j\ar[uu]^-{r_{X_0}}&K_0(A_3/A_1)\ar[uu]^-{r_{X_0}}\\
K^0(\Gamma_b^{tan}(X)|_{(X_0\times \{1\})^c})\ar[u]^-t\ar[r]_-{e_1^\partial}&K_0(A_2/A_0)\ar[r]_-j\ar[u]^-{q}& K_0(A_2/A_1)\ar[u]^-{q}\\
&&
}
\end{equation}
By commutativity of the diagram above
$$q(I)=J\in K_0(A_3/A_0),$$
if we denote by $I=e_1^\partial(\sigma^{nc}_{0})\in K_0(A_2/A_0)=K^0(\Gamma_b(X)|_{X_2\setminus X_0})$.

Notice also that $I\in Ker(K_0(A_2/A_1)\stackrel{\delta_2}{\longrightarrow}K_1(A_1/A_0))$ and that
the projection $q:K_0(A_2/A_0)\longrightarrow K_0(A_3/A_0)$ corresponds in conormal homology to the canonical projection $Ker\delta_2 \to \frac{Ker\delta_2}{Im\delta_3}$. We would then prove the theorem if we can prove that
\begin{equation}
j(I)=(Ind_{Fred,2}^{\overline{Y}}(\sigma^{nc}_{0,Y}))_{Y\in F_2^{cn}}\in K_0(A_2/A_1)\cong \bigoplus_{Y\in F_2}K^0(Y\times Y \times \mathbb{R}_2).
\end{equation}
Consider the following commutative diagram
\begin{equation}
\xymatrix{
K^0(\Gamma_b^{tan}(X)|_{(X_0\times \{1\})^c})\ar[d]_-a\ar[r]^-{e_1^\partial}&K_0(A_2/A_0)\ar[r]^-j\ar[d]^-j&K_0(A_2/A_1)\\
K^0(\Gamma_b^{tan}(X)|_{(X_1\times [0,1])\setminus X_0\times \{1\}})\ar[r]_-F&K_0(A_2/A_1)\ar[ru]_-=&
}
\end{equation}
where $F:K^0(\Gamma_b^{tan}(X)|_{(X_1\times [0,1])\setminus X_0\times \{1\}})\to K^0(\Gamma_b(X)|_{X_1\setminus X_0})=K_0(A_2/A_1)$
is the morphism induced by evaluation at $\epsilon=1$.
In particular we have that
\begin{equation}
j(I)=F(a(\sigma^{nc}_{0})\in K_0(A_2/A_1)\cong \bigoplus_{Y\in F_2}\mathbb{Z}\cdot Y
\end{equation}
We want then to show that the $Y-$component of $F(a(\sigma^{nc}_{0}))$, that we denote by $F(a(\sigma^{nc}_{0}))_Y$, is $Ind_{Fred,2}^{\overline{Y}}(\sigma^{nc}_{0,Y})$ for $Y\in F_2^{cn}(X)$ and zero otherwise. Now, the fact that $F(a(\sigma^{nc}_{0}))_Y=0$ for $Y\in F_2\setminus F_2^{cn}$ follows immediately from the fact that $F\circ a$ factors through $Ker\,\delta_2$ by the diagram above.

Let us fix then $Y\in F_2^{cn}(X)$, we have the following commutative diagram
\begin{equation}
\xymatrix{
K^0(\Gamma_b^{tan}(X)|_{(X_0\times \{1\})^c})\ar[r]^-a&K^0(\Gamma_b^{tan}(X)|_{(X_1\times [0,1])\setminus X_0\times \{1\}})\ar[d]_-{F}\ar[r]^-{r_{\overline{Y}}}&K^0(\Gamma_b^{tan}(X)|_{(\overline{Y}\times [0,1])\setminus \partial\overline{Y}\times \{1\}})\ar[d]^-{F_Y} \\
&K^0(\Gamma_b(X)|_{X_1\setminus X_0})\ar[r]_-{r_Y}&K^0(Y\times Y \times \mathbb{R}_2)
}
\end{equation}
where $r_{\overline{Y}}$ and $r_Y$ stand for the morphisms induced by the obvious restrictions indicated by the notations and where $F_Y$ is the morphism induced by evaluation at $\epsilon=1$.
Now, we have just to remark that 
$$K^0(\Gamma_b^{tan}(X)|_{(\overline{Y}\times [0,1])\setminus \partial\overline{Y}\times \{1\}})=K^0(\Gamma_b^{Fred}(\overline{Y})\times \mathbb{R}_2)$$
to obtain
\begin{equation}\label{Fredc}
F(a(\sigma^{nc}_{0}))_Y=Ind_{Fred,2}^{\overline{Y}}(\sigma^{nc}_{0,Y}) 
\end{equation}
and get the theorem statement.
\end{proof}

\subsubsection{Odd conormal index in codim 3}

Let $X$ be a mwc of codimension 3, we consider the odd conormal index morphism
\begin{equation}
\xymatrix{
K^1_{top}(^bT^*X)\ar[rr]^-{Ind_{cn,odd}^X}&&H_1^{pcn}(X)
}
\end{equation}
given as the composition of the odd analytic index 
\begin{equation}
\xymatrix{
K^1_{top}(^bT^*X)\ar[rr]^-{Ind_{ana,odd}^X}&&K_1(C^*(\Gamma_b(X)))
}
\end{equation}
and the isomorphism
\begin{equation}
K_1(C^*(\Gamma_b(X))\stackrel{T}{\longrightarrow} H_1^{pcn}(X)
\end{equation}
recalled above. Since 
\begin{equation}\label{decompositionH1}
H_1^{pcn}(X)=Ker\, \delta_3 \oplus H_1^{cn}(X)
\end{equation}
with $H_1^{cn}(X)=\frac{Ker\,\delta_1}{Im\,\delta_2}$, to compute explicitly $$Ind_{cn,odd}^X(\sigma)\in H_1^{pcn}(X)$$ 
for a given $\sigma\in K^1_{top}(^bT^*X)$, we will distinguish the two components
\begin{equation}
Ind_{cn,odd}^X(\sigma)_{\delta_3}\in Ker\, \delta_3
\end{equation}
and 
\begin{equation}
Ind_{cn,odd}^X(\sigma)_{\delta_1} \in H_1^{cn}(X)
\end{equation}
with respect to direct sum above.

\vspace{2mm}

{\bf The component in $Ker\,\delta_3$:} We need some terminology.

\begin{definition}
Let $Y\in F_3(X)$. Consider the canonical coefficient projection
$$C_3^{cn}(X)\stackrel{q_Y}{\longrightarrow}\mathbb{Z}$$
that computes the coefficient of $Y\otimes \epsilon_Y$. We say that $Y$ belongs to a corner cycle if there is $m\in Z_3^{cn}(X)$ such that $q_Y(m)\neq 0$.
We denote by $F_3^{cn}(X)\subset F_3(X)$ the subset of codimension 3 faces belonging to a cycle.
\end{definition}

Let $\sigma\in K^1_{top}(^bT^*X)$, for $Y\in F_3(X)$ there is a class given canonically by restriction
\begin{equation}
\sigma_Y\in K^1(T^*Y\times \mathbb{R}_3)
\end{equation}
and an associated $3$-suspended Atiyah-Singer index $Ind_{AS,3}^{Y}(\sigma_Y)\in \mathbb{Z}$.
We have the following computation in terms of corner cycles, whose proof follows the same lines as the proofs of the formulas for the odd conormal index for a codimension 1 mwc and for the even conormal index for a codimesion 2 mwc.

\begin{theorem}
Let  $\sigma \in K^1_{top}(^bT^*X)$, we have that
\begin{equation}
\sum_{Y\in F_3^{cn}(X)}Ind_{AS,odd}^Y(\sigma_Y)\cdot (Y\otimes \epsilon_Y)\in Ker\,\delta_3
\end{equation}
and
\begin{equation}
Ind_{cn,odd}^X(\sigma)_{\delta_3}=\sum_{Y\in F_3^{cn}(X)}Ind_{AS,3}^Y(\sigma_Y)\cdot (Y\otimes \epsilon_Y).
\end{equation}
\end{theorem}

\vspace{2mm}

{\bf The component in $H_1^{cn}(X)=\frac{Ker\,\delta_1}{Im\,\delta_2}$:} As for the cases of the odd conormal index in codim 2 and of the even conormal index in codim 3, the formula below will depend on some liftings, tamings, choices. Before stating the precise formula we need then to explain its main ingredients. 

Consider the following commutative diagram with exact columns
\begin{equation}
\xymatrix{
K^1(\Gamma_b(X)|_{X_1\times \{1\}})=K_1(A_3/A_1)\ar[r]^-{Id}_-=&K^1(\Gamma_b(X)|_{X_1\times \{1\}})=K_1(A_3/A_1)\\
&&\\
K^1(\Gamma_b^{tan}(X))\ar[r]^-{e_1}\ar[uu]^-{r_{X_1\times \{1\}}}&K_1(A_3)\ar[uu]^-{r_{X_1}}\\
K^1(\Gamma_b^{tan}(X)|_{(X_1\times \{1\})^c})\ar[u]^-t\ar[r]_-{e_1}&K_1(A_1).\ar[u]^-{q}\\
&&
}
\end{equation}

As recalled in the computation of $K_1(A_3)$ above, the morphism $r_{X_1}:K_1(A_3)\longrightarrow K_1(A_3/A_1)$ corresponds to $r_1: H^{pcn}_1(X)\to Ker\delta_3$, the projection onto the first factor of the direct sum decomposition (\ref{decompositionH1}). In particular $r_{X_1}$ is surjective, as well as $r_1:=r_{X_1\times \{1\}}$ by commutativity of the diagram above. Also, $K_1(A_3/A_1)\cong Ker\, \delta_3$ is free. We can hence consider a split injection $s_1:Ker\delta_3\to H^{pcn}_1(X)$ (the expression of the formula will depend on this choice), that is, a morphism $s_1$ with $r_1\circ s_1=Id_{Ker\,\delta_3}$.

Let  $\sigma \in K^1_{top}(^bT^*X)$ and $\sigma^{tan}\in K^1(\Gamma^{tan}_b)$ the unique element such that $e_0(\sigma^{tan})=\sigma$ under the isomorphism 
$$e_0:K^1(\Gamma^{tan}_b) \to K^1_{top}(^bT^*X)$$
induced by the evaluation at zero. Let 

\begin{equation}\label{definitionI1}
\sigma_1^{tan}=\sigma^{tan}-s_1(r_1(\sigma^{tan}))\in K^1(\Gamma^{tan}_b).
\end{equation}
By construction $r_1(\sigma_1^{tan})=0$ and hence we can then choose a lifting/taming $\sigma_1^{nc}\in K^1(\Gamma_b^{tan}(X)|_{(X_1\times \{1\})^c})$ such that $t(\sigma_1^{nc})=\sigma_1^{tan}$.

Let $Y\in F_1(X)$, in this case $\overline{Y}$ is a manifold with corners of codimension 2. Consider the following commutative diagram
\begin{equation}
\xymatrix{
K^1(\Gamma_b(X)|_{X_1})\ar[r]^-{r_{\partial \overline{Y}}}&K^1(\Gamma_b(\partial\overline{Y})\times \mathbb{R}_1)\\
&\\
K^1(\Gamma_b^{tan}(X))\ar[r]^-{r_{\overline{Y}}}\ar[uu]^-{r_{X_1\times \{1\}}}&K^1(\Gamma_b^{tan}(\overline{Y})\times\mathbb{R}_1)\ar[uu] _-{r_{\partial \overline{Y}\times\{1\}}}\\
K^1(\Gamma_b^{tan}(X)|_{(X_1\times \{1\})^c})\ar[u]^-t\ar[r]_-{r_{\overline{Y}}}&K^1(\Gamma_b^{Fred}(\overline{Y})\times \mathbb{R}_1)\ar[u]_-{\iota_0}
}
\end{equation}
where the horizontal morphisms are induced from the indicated restrictions to closed subgroupoids induced by restriction to $\overline{Y}$, the columns are exact induced from the indicated restrictions on the top.

We let $\sigma^{nc}_{1,Y}\in K^1(T_{nc}\overline{Y} \times \mathbb{R}_1)$ to be the image of $(r_Y)(\sigma^{nc}_1) \in K^1(\Gamma_b^{Fred}(\overline{Y})\times \mathbb{R}_1)$ under the isomorphism in K-theory
\begin{equation}
\xymatrix{
K^1(T_{nc}\overline{Y} \times \mathbb{R}_1)&K^1(\Gamma_b(\overline{Y})^{Fred}\times \mathbb{R}_1)\ar[l]^-{r_\partial}_-{\cong}.}
\end{equation}

We can hence consider the $1$-suspended Fredholm index 
$$Ind_{Fred,1}^{\overline{Y}}(\sigma^{nc}_{1,Y})\in \mathbb{Z}.$$ 
 
We can now state the computation of the odd conormal index in terms of these indices.

\begin{theorem}\label{Thmcnoddindexcodim3}
Let $\sigma\in K^1(^bT^*X)$, with the notations above we have that
\begin{equation}
\sum_{Y\in F_1(X)}Ind_{Fred,1}^{\overline{Y}}(\sigma^{nc}_{1,Y})\cdot (Y\otimes \epsilon_Y) \in Ker\, \delta_1
\end{equation}
and
\begin{equation}
Ind_{cn,odd}^X(\sigma)_{\delta_1}=\left[\sum_{Y\in F_1(X)}Ind_{Fred,1}^{\overline{Y}}(\sigma^{nc}_{1,Y})\cdot (Y\otimes \epsilon_Y)\right]
\end{equation}
where $[\cdot]$ denotes the class in $H_1^{cn}(X)$.
\end{theorem}
\begin{proof}
Recall that with the notation above we had $\sigma_1^{tan}\in K^1(\Gamma_b^{tan}(X))$ with $e_0(\sigma_1^{tan})=\sigma_1\in K^1(^bT^*X)$. Denote 
$$J=e_1(\sigma_1^{tan})\in K_1(A_3)\cong K_1(A_3/A_1)\oplus K_1(A_1)/Im(K_0(A_2/A_1)\to K_1(A_1)).$$
as we observed after (\ref{definitionI1}), we have by definition that $J \in K_1(A_1)/Im(K_0(A_2/A_1)\to K_1(A_1))$  which we identified above with $H_1^{cn}(X)=Ker\,\delta_1/Im\, \delta_2$.

Consider the following commutative diagram with exact columns
\begin{equation}
\xymatrix{
K^1(\Gamma_b(X)_{X_1})=K_1(A_3/A_1)\ar[r]^-{Id}_-=&K_1(A_3/A_1)\ar[r]^-{Id}_=&K_1(A_3/A_1)\\
&&\\
K^1(\Gamma_b^{tan}(X))\ar[r]^-{e_1}\ar[uu]^-{r_{X_1\times \{1\}}}&K_1(A_3)\ar[r]^-j\ar[uu]^-{r_{X_1}}&K_1(A_3/A_0)\ar[uu]^-{r_{X_1}}\\
K^1(\Gamma_b^{tan}(X)|_{(X_1\times \{1\})^c})\ar[u]^-t\ar[r]_-{e_1}&K_1(A_1)\ar[r]_-j\ar[u]^-{q}& K_1(A_1/A_0)\ar[u]^-{q}\\
&&
}
\end{equation}
By commutativity of the diagram above, we have
$$q(I)=J\in K_1(A_3),$$
if we denote by $I=e_1(\sigma_{1}^{nc})\in K_1(A_1)=K^1(\Gamma_b(X)|_{X_3\setminus X_1})$.

Notice also that $I\in Ker(K_1(A_1/A_0)\stackrel{\delta_1}{\longrightarrow}K_0(A_0))$ and that
the projection $q:K_1(A_1)\longrightarrow q(K_1(A_1))\subset K_1(A_2)$ corresponds in conormal homology to the canonical projection 
$$Ker\delta_1 \to \frac{Ker\delta_1}{Im\delta_2}.$$ 
We would then prove the theorem if we prove that
\begin{equation}
j(I)=(Ind_{Fred,1}^{\overline{Y}}(\sigma^{nc}_{1,Y}))_{Y\in F_1}\in K_1(A_1/A_0)\cong \bigoplus_{c\in F_1}K^1(Y\times Y \times \mathbb{R}).
\end{equation}
Consider the following commutative diagram
\begin{equation}
\xymatrix{
K^1(\Gamma_b^{tan}(X)|_{(X_1\times \{1\})^c})\ar[d]_-a\ar[r]^-{e_1}&K_1(A_1)\ar[r]^-j\ar[d]^-j&K_1(A_1/A_0)\\
K^1(\Gamma_b^{tan}(X)|_{(X_2\times [0,1])\setminus X_1\times \{1\}})\ar[r]_-F&K_1(A_1/A_0)\ar[ru]_-=&
}
\end{equation}
where $F:K^1(\Gamma_b^{tan}(X)|_{(X_2\times [0,1])\setminus X_1\times \{1\}})\to K^1(\Gamma_b(X)|_{X_2\setminus X_1})=K_1(A_1/A_0)$
is the morphism induced by evaluation at $\epsilon=1$.
In particular we have that
\begin{equation}
j(I)=F(a(\sigma_{1}^{nc}))\in K_1(A_1/A_0)\approx \bigoplus_{Y\in F_1}\mathbb{Z}\cdot Y
\end{equation}
We want then to show that the $Y-$component of $F(a(\sigma_1^{nc}))$, that we denote by $F(a(\sigma_1^{nc}))_Y$, is $Ind_{Fred,1}^{\overline{Y}}(\sigma^{nc}_{1,Y})$ for every $Y\in F_1(X)$.

Let us fix then $Y\in F_1(X)$, we have the following commutative diagram
\begin{equation}
\xymatrix{
K^1(\Gamma_b^{tan}(X)|_{(X_1\times \{1\})^c})\ar[r]^-a&K^1(\Gamma_b^{tan}(X)|_{(X_2\times [0,1])\setminus X_1\times \{1\}})\ar[d]_-{F}\ar[r]^-{r_{\overline{Y}}}&K^1(\Gamma_b^{tan}(X)|_{(\overline{Y}\times [0,1])\setminus \partial\overline{Y}\times \{1\}})\ar[d]^-{F_Y} \\
&K^1(\Gamma_b(X)|_{X_2\setminus X_1})\ar[r]_-{r_Y}&K^1(Y\times Y \times \mathbb{R}_1)
}
\end{equation}
where $r_{\overline{Y}}$ and $r_Y$ stand for the morphisms induced by the obvious restrictions indicated by the notations and where $F_Y$ is the morphism induced by evaluation at $\epsilon=1$.
We have just to remark that 
$$K^1(\Gamma_b^{tan}(X)|_{(\overline{Y}\times [0,1])\setminus \partial\overline{Y}\times \{1\}})=K^1(\Gamma_b^{Fred}(\overline{Y})\times \mathbb{R}_1)$$
to obtain
\begin{equation}\label{Fredc}
F(a(\sigma_1^{nc}))_Y=Ind_{Fred,1}^{\overline{Y}}(\sigma^{nc}_{1,Y})
\end{equation}
and get the theorem statement.
\end{proof}

\section{Fredholm anomalies in low codimensional closed manifold with corners}

\subsection{Ellipticity, fullellipticity and conormal obstructions for Fredholm boundary conditions}
\subsubsection{Ellipticity and Analytical Index morphisms}\label{subsectionAnalyticindexmorphism}

The analytical index morphism (of the manifold with embedded corners $X$)  takes it values in the group $K_0(\cK_b(X))$. It can be defined in two ways. First, we may consider  the connecting homomorphism $I$ of the exact sequence in $K$-theory associated with the short exact sequence of $C^*$-algebras:
\begin{equation}\label{bKses}
\xymatrix{
0\ar[r]&\cK_b(X)\ar[r]&\overline{\Psi_b^0(X)}\ar[r]^-{\sigma_b}&C(^bS^*X)\ar[r]&0.
}
\end{equation}
Then, if $[\sigma_b(D)]_1$ denotes the class in $K_1(C({}^bS^*X))$ of the principal symbol $\sigma_b(D)$ of an elliptic $b$-pseudodifferential $D$, we define the    analytical index   $\mathrm{Ind}_{\mathrm{an}}(D)$ of $D$ by 
\[
\mathrm{Ind}_{\mathrm{an}}(D)=I([\sigma_b(D)]_1)\in K_0(\cK_b(X)).
\]
Secondly, we can in a first step produce a $K_0$-class $[\sigma_b(D)]$ out of $\sigma_b(D)$: 
\begin{equation}
 [\sigma_b(D)] = \delta([\sigma_b(D)]_1)\in K_0(C_0({}^bT^*X))
\end{equation}
where $\delta$ is the connecting homomorphism of the exact sequence relating the vector and  sphere bundles:
\begin{equation}\label{bTS}
\xymatrix{
0\ar[r]& C_0({}^bT^*X)\ar[r]& C_0({}^bB^*X)\ar[r]&C({}^bS^*X)\ar[r]&0.
}
\end{equation}
Next,  we consider the exact sequence coming with the adiabatic deformation of $\Gamma_{b}(X)$:
\begin{equation}
\xymatrix{
0\ar[r]&C^*(\Gamma_b(X)\times (0,1])\ar[r]&C^*(\Gamma_b^{tan}(X))\ar[r]^-{r_0}&C^*(^bTX)\ar[r]&0,
}
\end{equation}
 in which the ideal is $K$-contractible. Using the shorthand notation $K^0_{top}(^bT^*X)$ for $K_0(C^*(^bTX))$, we set: 
\begin{equation}
Ind^a_X= r_1 \circ r_0^{-1} : K^0_{top}(^bT^*X)\longrightarrow K_0(\cK_b(X))
\end{equation}
where $r_1 : K_0(C^*(\Gamma_b^{tan}(X)))\to K_0(C^*(\Gamma_b(X))) $ is induced by the restriction morphism to $t=1$.
Applying a mapping cone argument to the exact sequence (\ref{bKses}) gives a commutative diagram
\begin{equation}
\xymatrix{
K_1(C(^bS^*X))\ar[rd]_-{\delta}\ar[rr]^-{I}&&K_0(\cK_b(X))\\
&K^0_{top}(^bT^*X)\ar[ru]_-{Ind^a_X}&
}
\end{equation}
Therefore we get, as announced:
\begin{equation}
  \mathrm{Ind}_{\mathrm{an}}(D) = Ind^a_X([\sigma_b(D)])
\end{equation}
 The map $Ind^a_X$ will be called  the {\it  Analytic Index morphism} of $X$.  A closely related homomorphism is the {\it  Boundary analytic Index morphism}, in which  the restriction to $X\times\{1\}$ is replaced  by the one to $\partial X\times\{1\}$, that is, we set:
 \begin{equation}
  Ind^\partial_X =  r_\partial \circ r_0^{-1}  : K_0(C_0(^bT^*X))\longrightarrow K_0(C^*(\Gamma_b(X)|_{\partial X})),
 \end{equation}
 where $r_\partial$ is induced by the   homomorphism  
 $C^*(\Gamma^{tan}_b(X))\longrightarrow C^*(\Gamma_b(X))|_{\partial X} $. We have of course 
 \begin{equation}
  Ind^\partial_X = r_{1,\partial}\circ Ind^a_X 
 \end{equation}
if $r_{1,\partial}$ denotes the map induced by the homomorphism  $C^*(\Gamma_b(X))\longrightarrow C^*(\Gamma_b(X)|_{\partial X})$. Since $r_{1,\partial}$ induces an isomorphism between $K_0$ groups (proposition 5.6 in \cite{CarLes}), both indices have the same meaning.

\subsubsection{Full ellipticity and the Fredholm Index morphism}\label{Fredsubsection}
To capture  the defect of Fredholmness of elliptic $b$-operators on $X$, we may introduce the algebra of full, or joint, symbols $\cA_\cF$ \cite{LMNpdo}.  If $F_1$ denotes the set of closed boundary hypersurfaces of $X$,  then the full symbol map is the $*$-homomorphism given by: 
\begin{equation}
 \sigma_F : \Psi^0(\Gamma_b(X))\ni P \longmapsto \Big( \sigma_b(P),(P\vert_H)_{H\in F_1}\Big) \in \cA_{\cF}.
\end{equation}
It gives rise to the exact sequence: 
\begin{equation}\label{Fredses}
\xymatrix{
0\ar[r]&\cK(X)\ar[r]&\overline{\Psi^0(\Gamma_b(X))}\ar[r]^-{\sigma_F}& \overline{\cA_{\cF}}\ar[r]&0
}
\end{equation}
where $\cK(X)$ is the algebra of compact operators on $L^2_b(X)$.
An operator $D\in \Psi^0(\Gamma_b(X))$ is said to be fully elliptic if $\sigma_F(D)$ is invertible.
 In \cite{Loya} (the statement also appears in \cite{MelPia}), it is proved that full ellipticity is equivalent to Fredholmness  on any $b$-Sobolev spaces $H^s_b(X)$.

For a given fully elliptic operator $D$, we denote by $\mathrm{Ind}_{\mathrm{Fred}}(D)$ its Fredholm index. We briefly recall how this integer is captured in $K$-theory. First, there is a natural isomorphism 
\begin{equation}
     K_0(\mu) \cong K_0(C^*(\mathcal{T}_{nc}X))
\end{equation}
between the $K$-theory of the obvious homomorphism $C(X)\stackrel{\mu}{\longrightarrow} \overline{\cA_{\cF}}$ and the $K$-theory of the noncommutative tangent space $\mathcal{T}_{nc}X$.  The former $K$-group captures stable homotopy classes of fully elliptic operators and the latter, which comes from deformation groupoid techniques,   classifies the noncommutative symbols $\sigma_{nc}(D)$ of  fully elliptic operators $D$.  

    Next, the same deformation  techniques give rise to a homomorphism:
 \begin{equation}\label{Fredmorph}
Ind^X_F : K^0(T_{nc}X)\longrightarrow
K_0(\cK(X))\simeq \mathbb{Z},
\end{equation}
which satisfies:
\begin{equation}
Ind_F^X([\sigma_{nc}(D)])= \mathrm{Ind}_{\mathrm{Fred}}(D),
\end{equation}
for any fully elliptic operator $D$. 

\subsubsection{Obstruction to full ellipticity and Fredholm perturbation property }\label{Obstructionsection}
 
 In order to analyse the obstruction to full ellipticity, we introduce Fredholm Perturbation Properties \cite{NisGauge}. 
\begin{definition} Let $D\in \Psi_b^m(X)$ be elliptic. We say that $D$ satisfies: 
\begin{itemize}
\item  the {\it Fredholm Perturbation Property} $(\cF\cP)$ if there is   $R\in \Psi_b^{-\infty}(X)$ such that $D+R$ is  fully elliptic. 
\item the  {\it stably  Fredholm Perturbation Property} $(\cS\cF\cP)$ if $D\oplus 1_H$ satisfies $(\cF\cP)$ for some identity operator $1_H$.
\item the {\it stably homotopic Fredholm Perturbation Property} $(\cH\cF\cP)$ if there is a fully elliptic operator $D'$ with $[\sigma_b(D')]=[\sigma_b(D)]\in K_0(C^*({}^bTX))$.
\end{itemize}
\end{definition}
We also say that $X$ satisfies the   {\it (resp. stably) Fredholm Perturbation Property}  if any elliptic $b$-operator on $X$ satisfies $(\cF\cP)$ (resp. $(\cS\cF\cP)$).

Property  $(\cF\cP)$ is  stronger than property $(\cS\cF\cP)$ which in turn is equivalent to property $(\cH\cF\cP)$ by  \cite[Proposition 4.3]{DebSkJGP}.  In  \cite{NSS2}, Nazaikinskii, Savin and Sternin characterized $(\cH\cF\cP)$
 for arbitrary manifolds with corners using an index map associated with their dual manifold construction. In \cite{CarLes} the result of \cite{NSS2} is rephrased in terms of deformation groupoids with the non trivial extra contribution of changing $(\cH\cF\cP)$ by $(\cS\cF\cP)$ thanks to   \cite[Proposition 4.3]{DebSkJGP}: 
  
\begin{theorem}\label{AnavsFredthm1}
Let $D$ be an elliptic $b$-pseudodifferential operator on a compact manifold with corners $X$. Then  $D$ satisfies $(\cS\cF\cP)$ if and only if 
$
 Ind_X^\partial([\sigma_b(D)])=0$ in  $K_0(C^*(\Gamma_b(X)|_{\partial X}))$.  \\ In particular, if $D$ satisfies $(\cF\cP)$ then its boundary analytic  index vanishes. 
\end{theorem}
This motivates the computation of $K_0(C^*(\Gamma_b(X)))\cong K_0(C^*(\Gamma_b(X)|_{\partial X}))$.

\subsection{Corner cycles and Fredholm anomalies in low codimension}

For $X$ a closed mwc of codimension lower or equal to three we can put the last theorem together with the isomorphism (computed in \cite{CarLes} and recalled above)

\begin{equation}
\xymatrix{
K^0(\Gamma_b(X))\ar[r]^-T_-\cong & H_{ev}^{cn}(X)
}
\end{equation}

to get the following corollary, which is altenatively a conormal homology version of the theorem above, 

\begin{corollary}
Let $D$ be an elliptic $b$-pseudodifferential operator on a compact manifold with corners $X$ of codimension lower or equal to three. Then  $D$ satisfies $(\cS\cF\cP)$ if and only if the even conormal index 
$
 Ind_X^{cn}([\sigma_b(D)])=0$ in  $H_{ev}^{cn}(X)$.  \\ In particular, if $D$ satisfies $(\cF\cP)$ then its even conormal index vanishes.
\end{corollary}

Now, using the explicit computations of $Ind_X^{cn}([\sigma_b(D)])$ for lower codimensions we can specialize the above corollary for codimensions 2 and 3 as follows:

\begin{corollary}[Fredholm anomalies in codimension 2]
Let $D$ be an elliptic $b$-pseudodifferential operator on a compact manifold with corners $X$ of codimension 2. Then  $D$ satisfies $(\cS\cF\cP)$ if and only if for every $Y\in F_2^{cn}(X)$ (face of codimension 2 
of $X$ belonging to a cycle, definition \ref{defF2cn}) the 2-suspended Atiyah-Singer index
\begin{equation}
Ind^Y_{AS,2}([\sigma_b(D)]_Y)=0.
\end{equation}
In particular, if $D$ satisfies $(\cF\cP)$ then all the 2-suspended Atiyah-Singer indices above vanish.
\end{corollary}

{\bf Terminology:} Remember we have been denoting by $$F_2^{cn}(X)=\{Y\in F_2(X):\exists m\in Z_2^{cn}(X) \text{ s.t. } q_Y(m)\neq 0\}.$$ We will now need another subset of the codimension 2 faces. Let 
$$F_2^{cn,\delta}(X)=\{Y\in F_2(X):\exists m\in Im(\delta_3) \text{ s.t. } q_Y(m)\neq 0\},$$
then of course $F_2^{cn,\delta}(X)\subset F_2^{cn}(X)$.

\begin{corollary}[Fredholm anomalies in codimension 3]
Let $D$ be an elliptic $b$-pseudodifferential operator on a compact manifold with corners $X$ of codimension 3. Then  $D$ satisfies $(\cS\cF\cP)$ if and only if 
\begin{enumerate}
\item For every $Y\in F_2^{cn}(X)\setminus F_2^{cn,\delta}(X)$ the 2-suspended Fredholm index 
\begin{equation}
Ind^{\overline{Y}}_{Fred,2}([\sigma_b(D)]_{0,Y}^{nc})=0,
\end{equation}
and
\item there is a 3-conormal chain $a\in C_3^{cn}(X)$ such that
\begin{equation}
\sum_{Y\in F_2^{cn,\delta}(X)}Ind_{Fred,2}^{\overline{Y}}([\sigma_b(D)]_{0,Y}^{nc})\cdot (Y\otimes \epsilon_Y)=\delta_3(a).
\end{equation}
\end{enumerate}
In particular, if $D$ satisfies $(\cF\cP)$ then the following two points above hold.
\end{corollary}

\subsection{Topological corner anomalies}

\subsubsection{Topological corner anomalies in codimension 2}

For $Y$ a smooth manifold we can define the topological Atiyah-Singer suspended index morphism as
\begin{equation}
Ind_{topAS,2}^{Y}:K^0(T^*Y\times \mathbb{R}_2)\to \mathbb{Q}
\end{equation}
defined by
\begin{equation}
Ind_{topAS,2}^{Y}(a):=\langle ch(\beta_2^{-1}(a)\wedge Td_Y, [Y]\rangle
\end{equation}
where $\beta_2:K^0(T^*Y)\to K^0(T^*Y\times \mathbb{R}_2)$ is the Bott isomorphism given by multiplication with the Bott element $\beta_2$ defined above. The classic Atiyah-Singer index theorem implies that
\begin{equation}
Ind_{AS,2}^{Y}=Ind_{topAS,2}^{Y},
\end{equation}
and in particular that $Ind_{topAS,2}^{Y}$ is integer valued. 
Now, for our purposes, the following is an immediate corollary of theorem \ref{cnindexcodim2}.

\begin{corollary}\label{topevencncodim2}
Let $X$ be a closed connected manifold with corners of codimension 2.
Let $P\in \overline{\Psi^0}(\Gamma_b(X))$ an elliptic operator with principal symbol class $[\sigma_P]\in K^0(^bT^*X)$, with the notations above  we have that
\begin{equation}
Ind_{cn,ev}^X([\sigma_P])=\sum_{Y\in F_2^{cn}(X)}Ind_{topAS,2}^{Y}([\sigma_P]_Y)\cdot (Y\otimes \epsilon_Y).
\end{equation}
\end{corollary}

As a corollary we obtain the following topological obstruction characterisation to the $(\cS\cF\cP)$-property for a given elliptic operator in codimension two:

\begin{theorem}[Topological anomalies in codimension 2]
Let $D$ be an elliptic $b$-pseudodifferential operator on a compact manifold with corners $X$ of codimension 2. Then  $D$ satisfies $(\cS\cF\cP)$ if and only if for every $Y\in F_2^{cn}(X)$ (face of codimension 2 
of $X$ belonging to a cycle, definition \ref{defF2cn}) the 2-suspended topological Atiyah-Singer index
\begin{equation}
\langle ch(\beta_2^{-1}([\sigma_b(D)]_Y))\wedge Td_Y, [Y]\rangle=0.
\end{equation}
In particular, if $D$ satisfies $(\cF\cP)$ then all the 2-suspended Atiyah-Singer indices above vanish.
\end{theorem}

\subsubsection{Topological corner anomalies in codimension 3}\label{subsectiontopanomalies3}
In the case of a mwc of codimension three the formula for the even conormal index morphism involves 2-suspended Fredholm indices of manifolds with boundary. We will briefly recall that for such manifolds these morphisms admit as well a topological computation.

Let $Z$ be a manifold with boundary $\partial Z$ defined by a global defining function $\rho$. Consider an embedding
$$i:Z \hookrightarrow \mathbb{R}^{N-1}$$
with $N$ even and an the associated clean embedding
$$j:Z\hookrightarrow \mathbb{R}^{N-1}\times \mathbb{R}$$
given by $j(z)=(i(z),\rho(z))$. Consider the smooth manifold
\begin{equation}
N_{sing}(j):=N(Z)\bigsqcup \mathbb{R}^{N-1}\times (0,1)
\end{equation}
obtained by gluying the manifolds with boundary $N(Z)$ and $D(\mathbb{R}^{N-1},\partial Z)|_{[0,1)}$ through their common boundary $N(\partial Z)$ (the gluying depends on the defininf function of the boundary). In \cite{CLM} we showed that there is Connes-Thom isomorphism
\begin{equation}
\xymatrix{
K^0(T_{nc}Z)\ar[r]^-{\mathcal{T}}_-\cong &K^0_{top}(N_{sing}(j)),
}
\end{equation}
then we defined
\begin{equation}\label{deftopfredindex}
Ind_{topFred}^Z:K^0(T_{nc}Z)\to \mathbb{Q}
\end{equation}
by 
\begin{equation}
Ind_{topFred}^Z(\sigma):=\int_{N_{sing(j)}}ch(\mathcal{T}(\sigma)).
\end{equation}
The main results in \cite{CLM}, theorem 4.4 and corollary 5.1, give an equality
\begin{equation}
Ind_{topFred}^Z=Ind_{Fred}^Z.
\end{equation} 
In particular, as for the AS case, $Ind_{topFred}^Z$ is integer valued and does not depend on the embedding. 
Up to a Bott isomorphism we get as well a 2-suspended topological Fredholm index $Ind_{topFred,2}^Z$ and the equality of suspended morphisms
\begin{equation}
Ind_{topFred,2}^Z=Ind_{Fred,2}^Z.
\end{equation}

Now, when applied to the computation of the even conormal morphism in codimension three we obtain:

\begin{corollary}\label{topevencncodim3}
Let $X$ be a closed connected manifold with corners of codimension 3.
Let $P\in \overline{\Psi^0}(\Gamma_b(X))$ an elliptic operator with principal symbol class $[\sigma_P]\in K^0(^bT^*X)$, with the notations above  we have that
\begin{equation}
Ind_{cn,ev}^X([\sigma_P])=\left[\sum_{Y\in F_2^{cn}(X)}Ind_{topFred,2}^{\overline{Y}}([\sigma_P]^{nc}_{0,Y})\cdot (Y\otimes \epsilon_Y)\right]
\end{equation}
where $[\cdot]$ denotes the class in $H_2^{cn}(X)$ and $[\sigma_P]^{nc}_{0,Y}\in K^0(T_{nc}\overline{Y}\times \mathbb{R}_2)$ are the noncommutative symbols defined in \ref{fixsigma}.
\end{corollary}

We get as an immediate corollary the theorem which states a topological formula for the obstructions to give Fredholm boundary conditions (up to stable property) on a given operator

\begin{theorem}[Topological corner anomalies in codimension 3]
Let $X$ be a closed connected manifold with corners of codimension 3. Then an elliptic $b$-operator $P\in \overline{\Psi^0}(\Gamma_b(X))$ has the $(SFP)$ if and only if 
\begin{enumerate}
\item For every $Y\in F_2^{cn}(X)\setminus F_2^{cn,\delta}(X)$ the 2-suspended topological Fredholm index 
\begin{equation}
Ind^{\overline{Y}}_{topFred,2}([\sigma_b(D)]_{0,Y}^{nc})=0,
\end{equation}
and
\item there is a 3-conormal chain $a\in C_3^{cn}(X)$ such that
\begin{equation}
\sum_{Y\in F_2^{cn,\delta}(X)}Ind_{topFred,2}^{\overline{Y}}([\sigma_b(D)]_{0,Y}^{nc})\cdot (Y\otimes \epsilon_Y)=\delta_3(a). 
\end{equation}
\end{enumerate}
In particular, if $P$ satisfies $(\cF\cP)$ then the following two points above hold.
\end{theorem}

\bibliographystyle{plain} 
\bibliography{CornersAindex} 

\end{document}